%% file: a.tex
\documentclass{article}

\usepackage{amsmath}
\usepackage{amssymb}
\usepackage{amsthm}
\usepackage{array}
\usepackage{bm}
\usepackage{color,soul}
\usepackage{comment}
\usepackage{footnote}
\usepackage{geometry}
\usepackage{graphicx}
\usepackage{hyperref}
\usepackage{cleveref}
\usepackage{makecell}
\usepackage[numbers]{natbib}
\usepackage{pgfplots}
\usepackage{ragged2e}
\usepackage{rotating}
\usepackage[labelformat=simple,listofformat=subsimple]{subfig}
\usepackage{tikz}

\makesavenoteenv{tabular}
\usetikzlibrary{decorations.pathreplacing}
\usetikzlibrary{backgrounds}
\usetikzlibrary{intersections}

\newtheorem{theorem}{Theorem}
\newtheorem{proposition}[theorem]{Proposition}%
\newtheorem{lemma}[theorem]{Lemma}
\newtheorem{corollary}[theorem]{Corollary}

\theoremstyle{remark}
\newtheorem{remark}{Remark}

\theoremstyle{definition}

\newtheorem{example}{Example}

\crefname{section}{Section}{Sections}
\crefname{subsection}{Section}{Sections}
\crefname{figure}{Fig.}{Figs.}
\crefname{subfigure}{Fig.}{Figs.}
\crefname{table}{Table}{Tables}
\crefname{lemma}{Lemma}{Lemmas}
\crefname{proposition}{Proposition}{Propositions}
\crefname{theorem}{Theorem}{Theorems}
\crefname{appendix}{Appendix}{Appendices}
\crefname{example}{Example}{Examples}

\newcommand{\email}[1]{\href{mailto:#1}{#1}}
\providecommand{\keywords}[1]
{
	\small	
	\textbf{\textit{Keywords---}} #1
}

\title{The numerical solution of fractional integral equations via orthogonal polynomials in fractional powers\thanks{The first author was supported by the Roth scholarship from the Department of Mathematics, Imperial College London. The second author was supported by the Leverhulme Trust Research Project Grant RPG-2019-144.}}
\author{Tianyi Pu\thanks{Department of Mathematics, Imperial College, London, United Kingdom (\email{tianyi.pu18@imperial.ac.uk}, \email{m.fasondini@imperial.ac.uk}).}\and Marco Fasondini\footnotemark[2]}
\date{June 2022}

\begin{document}

\maketitle

\begin{abstract}
	We present a spectral method for one-sided linear fractional integral equations on a closed interval that achieves exponentially fast convergence for a variety of equations, including ones with irrational order, multiple fractional orders, non-trivial variable coefficients, and initial-boundary conditions. The method uses an orthogonal basis that we refer to as  Jacobi fractional polynomials, which are obtained from an appropriate  change of variable in weighted classical Jacobi polynomials. New algorithms for building the matrices used to represent  fractional integration operators  are presented and compared. 
	Even though these algorithms  are unstable and require the use of high-precision computations, the spectral method nonetheless yields well-conditioned linear systems and is therefore stable and efficient. 
	For time-fractional heat and wave equations, we show that our method (which is not sparse but uses an orthogonal basis) outperforms a sparse spectral method (which uses a basis that is not orthogonal) due to its superior stability.
\end{abstract}
\keywords{fractional integral, spectral method, Jacobi polynomials, Riemann--Liouville, Caputo, Bagley--Torvik, Mittag--Leffler, high-precision}

\input{introduction}
\input{preliminaries}

\input{fjpolynomial}
\input{fioperator}
\input{spectral}
\input{conclusion}

\section*{Acknowledgements}
We thank Sheehan Olver for inspiring our interest in spectral methods for FIEs and FDEs and for many useful suggestions.

\section*{Code Availability}
\texttt{Julia} code for all figures presented in this paper are available at \url{https://github.com/putianyi889/JFP-demo}

\appendix
\input{ps}

\bibliography{reference}

\end{document}

%% file: introduction.tex
\section{Introduction}\label{intro}

Fractional differential equations (FDEs) and fractional integral equations (FIEs) have received a great deal of attention in the literature, not only because they appear in many scientific fields\footnote{For example, in physics \cite{hilfer2000applications,povstenko2015linear,benson2000application,bagley1983fractional}, chemistry \cite{mostafanejad2021fractional}, biology \cite{magin2010fractional}, finance \cite{FALLAHGOUL201712} and medical imaging \cite{li2015adaptive,treeby2010modeling}.}, but also because their solutions  are difficult to compute with traditional approaches. Numerous standard numerical methods\footnote{Examples include finite difference~\cite{cui2009compact}, finite element~\cite{ainsworth2018hybrid} and quadrature-based~\cite{diethelm1997algorithm} methods.} have been applied to FDEs and FIEs, however they only achieve algebraic convergence. This is because the fractional integral operators can introduce algebraic singularities that are not well approximated by polynomials. 

\begin{figure}[!htp]
    \centering
    \subfloat[Fractional integrals]{\includegraphics[width=0.4\textwidth]{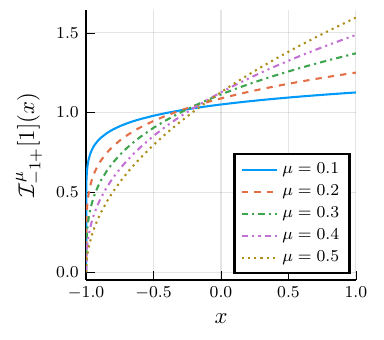}\label{flo51}}
    \quad
    \subfloat[Errors]{\includegraphics[width=0.4\textwidth]{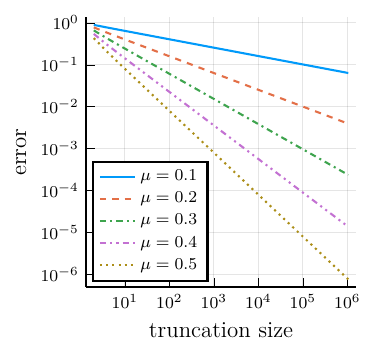}\label{flo52}}
    \caption{Algebraic convergence of polynomial approximants to a fractional integral. \protect\subref{flo51}: Plots of fractional integrals of the unit constant function, $\mathcal{I}^{\mu}_{-1^+}[1](x)$.
    \protect\subref{flo52}: Errors from Chebyshev polynomial approximations to $\mathcal{I}^{\mu}_{-1^+}[1](x)$, which converge as $\mathcal{O}(n^{-2\mu})$, where $n+1$ is the truncation size.}\label{flo53}
\end{figure}

The fractional integral of order $\mu > 0$ is a singular integral defined as
\begin{equation*}
    \mathcal{I}^\mu_{a^+}\varphi(x):=\frac{1}{\Gamma(\mu)}\int_{a}^{x}\frac{\varphi(t)}{(x-t)^{1-\mu}}dt\quad\text{ or }\quad \mathcal{I}^\mu_{b^-}\varphi(x):=\frac{1}{\Gamma(\mu)}\int_{x}^{b}\frac{\varphi(t)}{(t-x)^{1-\mu}}dt,
\end{equation*}
with $a < x < b$ and where $\mathcal{I}^\mu_{a^+}$ and $\mathcal{I}^\mu_{b^-}$ are referred to as, respectively, the left-sided and right-sided fractional integrals. 
Consider the left-sided fractional integral of the unit constant function, $\mathcal{I}^{\mu}_{-1^+}[1](x)$, which is $(1 + x)^{\mu}/\Gamma(1+\mu)$ (see \eqref{eq4}). Hence, $\mathcal{I}^{\mu}_{-1^+}[1](x)$ has an algebraic singularity at $x = -1$ as illustrated by \cref{flo51} and the Chebyshev polynomial approximation of $\mathcal{I}^{\mu}_{-1^+}[1](x)$ on $[-1,1]$ converges at the slow algebraic rate of $\mathcal{O}(n^{-2\mu})$ as indicated by \cref{flo52}, where $n$ is the polynomial degree. By contrast, with the Jacobi fractional polynomial (JFP) basis that we shall introduce, algebraic singularities are incorporated into the basis, which allow us to compute fractional integrals of weighted polynomials in fractional powers exactly.  

Fractional derivatives have different definitions, among which the Riemann--Liouville and Caputo types are commonly used. We denote them by, respectively, $\mathcal{D}^\nu_{\mathrm{RL}}$ and $\mathcal{D}^\nu_{\mathrm{C}}$, where $\nu>0$ and they are defined as compositions of fractional integral and standard derivative operators (but in different orders) as follows~\cite[(2.22)]{kilbas1993fractional},
\begin{equation} 
\mathcal{D}^\nu_{\mathrm{RL}}\varphi(x):=\mathcal{D}^m\mathcal{I}^{\mu}\varphi(x),\quad \mathcal{D}^\nu_{\mathrm{C}}\varphi(x):=\mathcal{I}^{\mu}\left[\mathcal{D}^m\varphi\right](x),\quad m =  \lceil\nu\rceil,\quad \mu=m-\nu.  \label{eq:RLCdefs}
\end{equation}
Here $\mathcal{D}^m$ denotes the $m$-th order (standard) differential operator: $\mathcal{D}^m\varphi(x)=\frac{d^m}{dx^m}\varphi(x)$, and $\mathcal{I}^{\mu}$ can be either a left-sided or right-sided fractional integral.
In this paper we focus on problems involving fractional integral operators. We do consider FDEs, however we reformulate them as FIEs, similar to the integral reformulation method discussed in~\cite[Section 5]{olver2020fast} in which ODEs are reformulated as ordinary integral equations.
In future work we shall construct differentiation matrices for JFP bases, which will make our method, which we refer to as the \emph{JFP method}, directly applicable to FDEs, without the need for integral reformulation.

We now give a brief overview of existing methods that can achieve exponential convergence for FDEs and FIEs and compare them to the method we shall present. All of these approaches incorporate algebraic singularities into the basis functions.

The fractional collocation methods in \cite{jiao2016well,zayernouri2014fractional} and the fractional Galerkin method in \cite{mao2016efficient}  achieve spectral convergence when the solution is smooth. However, in addition to giving rise to dense linear systems, these methods revert back to algebraic convergence if the solution has a singularity at the left endpoint, which is expected in practice (e.g., \cref{exa2}). We shall overcome this problem in the JFP method by including an extra singularity into the basis such that it captures the singularity of the solution. This idea is also applicable to the methods in \cite{jiao2016well,zayernouri2014fractional,mao2016efficient} and should restore spectral convergence to the examples where they achieved algebraic convergence.

Our work is strongly influenced by the method proposed by Hale and Olver~\cite{hale2018fast}, in which direct sums of appropriately weighted Jacobi polynomials are used as basis functions. We shall refer to this method as the \emph{sum space method} and we note that the basis functions are related to the ``generalized Jacobi functions'' of~\cite{chen2016generalized} and the ``polyfractonomials'' of~\cite{zayernouri2014fractional}. Similar to the ultraspherical spectral method~\cite{olver2013fast} for ODEs, in the sum space method, the domain and range of operators in FDEs and FIEs are represented in different  bases to ensure that the matrix representations of operators are banded. This enables a fast algorithm with linear complexity and exponential convergence for a wide range of problems. However, we shall find in \cref{spec} that the sum space method leads to tremendously large expansion coefficients (larger than $10^{100}$)  for the solution of an important family of FIEs that arises in the solution to the time-fractional heat/wave equation ($\mathcal{D}^{\mu}_t u = \mathcal{D}^2_x u$). These large sum space coefficients require very high precision for the accurate computation of solutions and is thus very expensive.  By contrast, this family of FIEs poses no difficulties to the JFP method whose solution coefficients are bounded below $1$. We shall clarify this difference in the performance of the two methods by using bounds on Jacobi expansion coefficients for analytic functions to bound the JFP coefficients and by illustrating that the largest expansion coefficient in the sum space basis grows at the same (super-exponential) rate as the largest coefficient of the power series expansion of the solution. 

We consider our method to be the successor of that of Bhrawy and Zaky~\cite{bhrawy2016shifted}. They applied a change of variables to classical Jacobi polynomials such that the algebraic singularities of the resulting basis, the JFP basis (which is called thus for reasons we explain in \cref{jfp}), conform to those of the solution\footnote{The method of  Bhrawy and Zaky can be considered to be an extension of that of Kazem, et al.~\cite{kazem2013fractional}, which used Legendre polynomials in fractional powers as basis functions.}. The JFP basis inherits many desirable properties of classical Jacobi polynomials, including orthogonality, see~\cite{bhrawy2016shifted}. Therefore, it can potentially be used to solve FIEs and FDEs just as classical orthogonal polynomials are used to solve ODEs. We believe this approach has not been sufficiently analyzed nor developed to its full potential. In this paper, we develop new algorithms for computing the entries of matrices that represent the action of fractional integration operators on the JFP basis. We shall refer to these matrices as fractional integration matrices. We also illustrate that the algorithms for fractional integration matrices (including the one used by Bhrawy and Zaky) are unstable but can be ``pseudo-stabilized" by integrating high-precision computing~\cite{bailey2012high} automatically such that the algorithms output accurate results.  We also investigate empirically the performance of the parameter-dependent JFP method for a wider range of parameters than those considered in~\cite{bhrawy2016shifted}. 

Our pseudo-stabilization technique discussed  in \cref{ps} is essential to the scalability of the JFP method and thus also for its application to practical computational problems. We emphasize that high-precision computations are required only for the computation of fractional integration matrices and not for the solution of the resulting linear systems. In addition, we do not require any of the inputs to the FIE or FDE (e.g., variable coefficients, the function on the right-hand side of the equation, boundary conditions, etc.) to be  computable to high-precision accuracy.

The following is an outline of the paper: We introduce the basic constituents of the JFP method in \cref{sect:prelims,jfp} (matrix representations of operators on quasimatrices, high-precision floating-point numbers, the JFP basis, etc.). We then focus on the properties (\cref{fio}) and computation (\cref{sect:fioalgs}) of fractional integration operators and matrices acting on the JFP basis. The JFP method is used in \cref{spec} to solve a variety of FIEs\footnote{\texttt{Julia} code for the examples in \cref{spec} are available at \url{https://github.com/putianyi889/JFP-demo}}, including an FDE and fractional PDE reformulated as FIEs, and comparisons are made with the sum space method. We conclude the paper with a summary and a discussion of topics for future work. \cref{ps} is devoted to the above-mentioned pseudo-stabilization of an algorithm for computing fractional integration matrices. 

%% file: preliminaries.tex
\section{Preliminaries}\label{sect:prelims}
Henceforth, we consider left-handed fractional integrals and derivatives on compact intervals, which we can take to be defined on $[-1,1]$, without loss of generality. To simplify the notation, we let
\begin{equation}\label{eq1}
    \mathcal{I}^\mu\varphi(x):=\frac{1}{\Gamma(\mu)}\int_{-1}^x\frac{\varphi(t)}{(x-t)^{1-\mu}}dt
\end{equation}
and let $\mathcal{I}$ denote $\mathcal{I}^1$, the standard integral operator. To avoid ambiguity, we shall denote (infinite) identity matrices by $\mathit{1}$.
We shall make use of the fact that fractional integral operators form a semigroup~\cite[(2.21)]{kilbas1993fractional}:
\begin{equation}\label{eq20}
    \mathcal{I}^\mu\mathcal{I}^\nu=\mathcal{I}^\nu\mathcal{I}^\mu=\mathcal{I}^{\mu+\nu}.
\end{equation}

\input{pre_qn}
\input{pre_jp}
\input{pre_hp}

%% file: pre_qn.tex
\subsection{Quasimatrix notation and fractional monomial bases} As in~\cite{hale2018fast,olver2020fast}, we shall adopt quasimatrix notation because of its convenience. 
A quasimatrix can be thought of as a matrix with a function in every column rather than a vector. For instance, we let
\begin{equation*}\mathbf{M}_0(x) =
\begingroup 
\setlength\arraycolsep{8pt}
 \left(\begin{array}{c | c | c | c | c}
 1 & (1+x) & (1+x)^2 & (1+x)^3 & \cdots
 \end{array}
\right)
\endgroup,
\end{equation*}
where $x \in [-1, 1]$; $\mathbf{M}_0(x)$ (or simply $\mathbf{M}_0$) is an example of a quasimatrix that has the (shifted) monomial basis functions in its columns, ordered by degrees, hence it has a countable infinity of columns and we write $\mathbf{M}_0 \in \mathbb{R}
^{[-1, 1] \times \mathbb{N}_0}$ (similar to the way in which we might state $A \in \mathbb{R}^{m\times n}$ for a matrix $A$). Let $\gamma \in \mathbb{R}$, then we define
\begin{equation*}
\mathbf{M}_0^{\gamma}(x) = 
\begingroup 
\setlength\arraycolsep{8pt}
\left(
\begin{array}{c | c | c | c | c}
1 & (1 + x)^{\gamma} & (1 + x)^{2\gamma} & (1 + x)^{3\gamma} & \cdots
\end{array}
\right).
\endgroup
\end{equation*}
The columns of $\mathbf{M}_0^{\gamma}$ can be viewed as elements of a fractional monomial basis if $\gamma$ has a positive non-integer value. We shall also make use of a normalized version of $\mathbf{M}^{\gamma}_{0}$, obtained by replacing $(1+x) \mapsto (1+x)/2$, i.e.,
\begin{equation}
\widetilde{\mathbf{M}}^{\gamma}_{0}(x) = 
\begingroup 
\setlength\arraycolsep{8pt}
\left(
\begin{array}{c | c | c |  c}
1 & \left(\frac{1 + x}{2}\right)^{\gamma} & \left(\frac{1 + x}{2}\right)^{2\gamma}  & \cdots
\end{array}
\right). \endgroup
 \label{eq:normM0}
\end{equation} 
 We denote the quasimatrix $\mathbf{M}_0^{\gamma}$ weighted by $(1 + x)^{\delta}$, i.e., $(1 + x)^{\delta}\mathbf{M}_0^{\gamma}(x)$,  by $\mathbf{M}^{\gamma}_{\delta}(x)$, hence
 \begin{equation}
 \begingroup 
\setlength\arraycolsep{5pt}
\mathbf{M}^{\gamma}_{\delta}(x) = (1 + x)^{\delta}\mathbf{M}_0^{\gamma}(x) = \left(
\begin{array}{c | c | c |  c}
(1 + x)^{\delta} & (1 + x)^{\gamma + \delta} & (1 + x)^{2\gamma + \delta}  & \cdots
\end{array}
\right).\endgroup
  \label{eq:Mgd}
\end{equation}

We can express the action of the fractional integral operator on $\mathbf{M}^{\gamma}_{\delta}$ using the fact that~\cite[(2.44)]{kilbas1993fractional}:
\begin{equation}\label{eq4}
    \mathcal{I}^\mu[(1+\diamond)^\beta](x)=\frac{\Gamma(\beta+1)}{\Gamma(\beta+\mu+1)}(1+x)^{\beta+\mu}
\end{equation}
where $\beta > -1$, to ensure the integral exists. Throughout, we shall use $\diamond$ as a dummy variable. The following result is an immediate corollary of \eqref{eq4}:
\begin{proposition}\label{prop:fiomono}
Let $\gamma > 0$ and $\mu = k\gamma$, where $k$ is a positive integer (i.e., $k \in \mathbb{N}_+$), 
 then $\mathcal{I}^{\mu}$ maps $\mathbf{M}^{\gamma}_{\delta}$ to itself via an infinite matrix whose $k$-th subdiagonal is nonzero, i.e.,
$\mathcal{I}^{\mu} \mathbf{M}^{\gamma}_{\delta} = \mathbf{M}^{\gamma}_{\delta}\Lambda^k_{\delta,\gamma}$,
where 
\begin{equation}
\Lambda^k_{\delta,\gamma} = \left(
\begin{array}{l}
\left.
\begin{array}{c}
0  \\
\vdots \\
0 
\end{array}\right\}k \text{ zeros} \\
\\
\hspace*{-0.3 cm}\begin{array}{l l l l}
\frac{\Gamma(\delta + 1)}{\Gamma(\delta+\mu+1)} & & & \\
 &  \frac{\Gamma(\delta+\gamma + 1)}{\Gamma(\delta+\gamma+\mu+1)} & & \\
 & & \frac{\Gamma(\delta+2\gamma + 1)}{\Gamma( \delta+2\gamma+\mu+1)} & \\
  & & & \ddots
\end{array}
\end{array}
\right), \qquad \Lambda_{\delta,\gamma}=\Lambda^1_{\delta,\gamma}. \label{eq:Lkdiag}
\end{equation}
\end{proposition}
We only consider the case $\mu = k\gamma$ because, as we shall find in \cref{spec}, this is a necessary condition for the solutions to the FIEs to have an expansion in the JFP basis.

We say that a matrix $A$ has bandwidths $(\lambda, \omega)$ (or $\mathrm{bandwidths}(A) = (\lambda, \omega)$) if the entries of $A$ satisfy $A_{i,j} = 0$ for  $i - j > \lambda$ and $j - i > \omega$. For example, if $A$ is upper triangular, then we say that $\mathrm{bandwidths}(A) = (0, \infty)$ if $A$ is infinite and $\mathrm{bandwidths}(A) = (0, N)$ if $A$ is $N \times N$. As another example, regardless of whether the matrix \eqref{eq:Lkdiag} is infinite or $N \times N$, it  has bandwidths $(k,-k)$ .

For quasimatrices $\mathbf{P}$ and $\mathbf{Q}$ and a weight function $w$, we let 
$\langle \mathbf{P}, \mathbf{Q} \rangle_{w}$
denote a matrix of $L^2(w)$ inner products between the functions forming the columns of  $\mathbf{P}$ and $\mathbf{Q}$. As an example we shall use later, let $\mathbf{P} = \mathbf{Q} = \mathbf{M}_0$ and  let $w$ be the Jacobi weight function, $w = w_{\alpha, \beta}$, where
\begin{equation}
w_{\alpha, \beta}(x) = (1-x)^{\alpha}(1+x)^{\beta}, \qquad x \in [-1, 1], \qquad \alpha, \beta > -1.  \label{eq:jacw}
\end{equation}
Also, let
\begin{equation}
\langle f,g\rangle_{w_{\alpha,\beta}}:=\int_{-1}^1f(x)g(x)w_{\alpha,\beta}(x)dx,  \label{eq:jacip}
\end{equation}
and define $\bm{e}_k$, $k \geq 0$ to be the basis vector
\begin{equation*}
\bm{e}_k = \Big( \underbrace{0 \: \cdots \: 0}_{k \text{ zeros}} \:\: 1  \:\: 0  \:\: 0 \:\: \cdots \Big)^{\top},
\end{equation*}
 then $\mathbf{M}_0\bm{e}_k$ is the function in column $k$ of $\mathbf{M}_0$ (i.e., $(1+x)^k$) and  the $(i,j)$-entry (with $i, j \geq 0$) of the infinite matrix $\langle \mathbf{M}_0, \mathbf{M}_0\rangle_{w_{\alpha,\beta}}$, which we denote by $M^{(\alpha,\beta)}$, is
  \begin{align}
M^{(\alpha,\beta)}_{i,j} = \left(\langle \mathbf{M}_0, \mathbf{M}_0\rangle_{w_{\alpha,\beta}} \right)_{i,j} &= \int_{-1}^{1} \left(\mathbf{M}_0^{\top} \mathbf{M}_0 \right)_{i,j} w_{\alpha,\beta}\, d x \label{eq:Mab}  \\
 & = \langle \mathbf{M}_0\bm{e}_i, \mathbf{M}_0\bm{e}_j \rangle_{w_{\alpha,\beta}} \label{eq:Mabent}	 \\
& =  \int_{-1}^{1}(1+x)^{\beta + i + j}(1 - x)^{\alpha} d x \notag \\
	& = 2^{\alpha+\beta+i+j+1}\mathrm{B}(\beta+i+j+1,\alpha+1), \notag
\end{align}
where $\mathrm{B}$ is the beta function. Since $\mathbf{M}_0 \in \mathbb{R}^{[-1, 1] \times \mathbb{N}_0}$, in \eqref{eq:Mab} we interpret  $\mathbf{M}_0^{\top} \mathbf{M}_0 \in \mathbb{R}^{ \mathbb{N}_0 \times \mathbb{N}_0}$ as the outer product of two vectors\footnote{Not to be confused with the notation in~\cite{olver2020fast}, where $\mathbf{P}^{\top} w \mathbf{Q}$ is used to denote  $\langle \mathbf{P}, \mathbf{Q} \rangle_{w}$. } whose entries are functions on $[-1, 1]$. In the general case, the $(i,j)$ entry  of $\langle \mathbf{P}, \mathbf{Q} \rangle_{w}$ is defined as $ \langle \mathbf{P}\bm{e}_i, \mathbf{Q}\bm{e}_j \rangle_{w}$ (cf.\eqref{eq:Mabent}).  If $\mathbf{P} = \mathbf{Q}$, then the (infinite and symmetric) matrix $\langle \mathbf{P}, \mathbf{P} \rangle_{w}$ is known as a Gram matrix~\cite[p.~441]{horn2012matrix}, hence $M^{(\alpha,\beta)}$ is a Gram matrix.

%% file: pre_jp.tex
\subsection{Jacobi polynomials}\label{jp}

We let the columns of the quasimatrix $\mathbf{P}^{(\alpha,\beta)}$ consist of Jacobi polynomials, ordered by degrees,
\begin{equation*}
\mathbf{P}^{(\alpha,\beta)}(x) =\begingroup 
\setlength\arraycolsep{8pt} \left(
\begin{array}{c | c | c}
P_0^{(\alpha,\beta)}(x)& P_1^{(\alpha,\beta)}(x) & \cdots
\end{array}
\right)\endgroup,
\end{equation*}
where the  $P_n^{(\alpha,\beta)}$ are orthogonal with respect to the inner product \eqref{eq:jacip} with\footnote{In~\cite[4.22]{szeg1939orthogonal} and \cite{alfaro2002orthogonality, chen2016generalized} Jacobi polynomials are considered for $\alpha \leq -1$ or $\beta \leq -1$.} $\alpha, \beta > -1$. The $P_n^{(\alpha,\beta)}$ are related to the (shifted) monomial basis functions according to
 \begin{equation*}
P^{(\alpha,\beta)}_n(x)
=\sum_{k=0}^n\frac{(k+\beta+1)_{n-k}(n+\alpha+\beta+1)_k}{(n-k)!k!}\left(\frac{x+1}{2}\right)^k(-1)^{n-k}, 
\end{equation*}
which follows from the identities \cite[(18.5.7)]{NIST:DLMF} and $P^{(\beta,\alpha)}_n(x) = (-1)^nP^{(\beta,\alpha)}_n(-x)$ (see \cite[Table 18.6.1]{NIST:DLMF}). Hence, 
\begin{equation}\label{eq6}
    \mathbf{P}^{(\alpha,\beta)}=\mathbf{M}_0 C^{(\alpha,\beta)},
\end{equation}
where $C^{(\alpha,\beta)}$ is an upper triangular matrix with entries 
\begin{equation}
C^{(\alpha,\beta)}_{k,n} =\frac{(-1)^{n-k}(k+\beta+1)_{n-k}(n+\alpha+\beta+1)_k}{2^k(n-k)!k!}, \qquad 0 \leq k \leq n.  \label{eq:Ckn}
\end{equation}

We can use the orthogonality of the Jacobi polynomials to obtain the inverse of $C^{(\alpha,\beta)}$. Define $\left\| \mathbf{P}^{(\alpha, \beta)} \right\|^2 = \langle \mathbf{P}^{(\alpha, \beta)}, \mathbf{P}^{(\alpha, \beta)}\rangle_{  w_{\alpha,\beta} }$,
 then $\left\| \mathbf{P}^{(\alpha, \beta)} \right\|^2$ is a diagonal matrix with the squared norms of the Jacobi polynomials on the main diagonal, which are given in~\cite[Table 18.3.1]{NIST:DLMF}.  Using \labelcref{eq6,eq:Mab}, it follows that
 \begin{equation}
\left\| \mathbf{P}^{(\alpha, \beta)} \right\|^2 = \left(C^{(\alpha,\beta)}\right)^{\top}M^{(\alpha,\beta)}C^{(\alpha,\beta)},     \label{eq:Pabnorm}
\end{equation}
 and thus
\begin{equation}
\left(C^{(\alpha,\beta)}\right)^{-1} =\left(   \left\| \mathbf{P}^{(\alpha, \beta)} \right\|^2 \right)^{-1} \left(C^{(\alpha,\beta)}\right)^{\top}M^{(\alpha,\beta)}.  \label{eq:Cabinvdef}
\end{equation}

Since $C^{(\alpha,\beta)}$ is upper triangular, so is $\left(C^{(\alpha,\beta)}\right)^{-1}$ and since the entries of $M^{(\alpha,\beta)}$ are known (see \eqref{eq:Mab}), the $(i,j)$ entry (with $i, j \geq 0$) of $\left(C^{(\alpha,\beta)}\right)^{-1}$ can be computed as a finite sum of $i+1$ terms\footnote{\eqref{eq:Pabnorm} can be expressed as  $M^{(\alpha,\beta)} = \left(C^{(\alpha,\beta)}\right)^{-\top} \left\| \mathbf{P}^{(\alpha, \beta)} \right\|^2 \left(C^{(\alpha,\beta)}\right)^{-1}$, which is an LDL factorisation of the symmetric positive-definite matrix $M^{(\alpha,\beta)}$.}.

We shall use the matrices acting on Jacobi polynomials given in \cref{flo1} as building blocks to construct matrix representations of the operators we need. For example, the integration matrix (which is a  representation of the standard  integration operator) is diagonal and follows from the weighted differentiation matrix in \cref{flo1}:
\begin{equation}
    \mathcal{I}[(1+\diamond)^\beta\mathbf{P}^{(\alpha,\beta)}](x)=(1+x)^{\beta+1}\mathbf{P}^{(\alpha-1,\beta+1)}(x)\left(W^{(\alpha,\beta)}_{(\alpha-1,\beta+1)}\right)^{-1}.  \label{eq:wintop}
\end{equation}
We shall also use the following matrix representation of the multiplication operator (which we shall refer to as a multiplication matrix):
\begin{proposition}\label{thm2} We have that
\begin{equation*}
(1+x)\mathbf{P}^{(\alpha,\beta)}(x)=\mathbf{P}^{(\alpha,\beta)}(x)L_{(\alpha,\beta+1)}^{(\alpha,\beta)}R_{(\alpha,\beta)}^{(\alpha,\beta+1)},
\end{equation*}
where the multiplication matrix $L_{(\alpha,\beta+1)}^{(\alpha,\beta)}R_{(\alpha,\beta)}^{(\alpha,\beta+1)}$ has bandwidths $(1,1)$.
\end{proposition}
\begin{proof}
By the actions of the matrices in \cref{flo1},
\begin{align*}
(1+x)\mathbf{P}^{(\alpha,\beta)}(x)  = (1+x)\mathbf{P}^{(\alpha,\beta+1)}(x)R_{(\alpha,\beta)}^{(\alpha,\beta+1)} =  \mathbf{P}^{(\alpha,\beta)}(x)L_{(\alpha,\beta+1)}^{(\alpha,\beta)}R_{(\alpha,\beta)}^{(\alpha,\beta+1)},
\end{align*}
and $L_{(\alpha,\beta+1)}^{(\alpha,\beta)}R_{(\alpha,\beta)}^{(\alpha,\beta+1)}$ has bandwidths $(1, 1)$ since it is a product of matrices with bandwidths $(1, 0)$ and $(0, 1)$.
\end{proof}

\begin{table}
\begin{tabular}{ m{0.25\textwidth} | m{0.48\textwidth} | m{0.15\textwidth}}
\RaggedRight
\textbf{Operator}  
            & \textbf{Action via matrix representation}
                & \textbf{Bandwidths} \\\hline
\footnotesize{Conversion} 
            &\footnotesize{$\mathbf{P}^{(\alpha,\beta)}=\mathbf{P}^{(\alpha+k,\beta+j)}R_{(\alpha,\beta)}^{(\alpha+k,\beta+j)}$}
                & \footnotesize{$(0,k+j)$}
                \\\hline
\footnotesize{Weighted conversion}
            &\footnotesize{ $w_{k,j}\mathbf{P}^{(\alpha+k,\beta+j)}=\mathbf{P}^{(\alpha,\beta)}L_{(\alpha+k,\beta+j)}^{(\alpha,\beta)}$}
                & \footnotesize{$(k+j,0)$}\\\hline
\footnotesize{Weighted differentiation}
           &\footnotesize{$\mathcal{D}[w_{0,\beta+1}\mathbf{P}^{(\alpha,\beta+1)}]=w_{0,\beta}\mathbf{P}^{(\alpha+1,\beta)}W_{(\alpha,\beta+1)}^{(\alpha+1,\beta)}$}
                & \footnotesize{Diagonal}
                \\\hline
\end{tabular}
\caption{The action of operators on quasimatrices of Jacobi polynomials via their banded matrix representations. In the final column, $k, j \in \mathbb{N}_0$. See~\cite[Appendix A]{olver2020fast} for the entries of these matrices.}\label{flo1}
\end{table}

%% file: pre_hp.tex
\subsection{High-precision computation}\label{sect:hp}

We shall use high-precision computation~\cite{bailey2012high} in the \texttt{Julia} programming language\footnote{See \url{https://julialang.org/} and \url{https://docs.julialang.org/en/v1/manual/integers-and-floating-point-numbers/\#Arbitrary-Precision-Arithmetic}} due to the instability of our algorithms for computing fractional integration matrices, which are discussed in \cref{sect:fioalgs}. As in~\cite{10.1145/1236463.1236468}, ``$q$-bit precision'' means the machine epsilon (or relative accuracy of the floating-point numbers we compute with) is $2^{1-q}$. For example, the widely-used double-precision arithmetic has 53-bit precision and a machine epsilon of 2.22e-16. Some of the commonly used precisions are given in \cref{flo36}.
\begin{table}[!htp]
    \centering
    \begin{tabular}{|c|c|}
        \hline
         \textbf{Precision}&\textbf{Machine epsilon}\\\hline
         \footnotesize{double precision} & \footnotesize{2.22e-16}\\
         \footnotesize{128-bit precision} & \footnotesize{5.88e-39}\\
         \footnotesize{256-bit precision} & \footnotesize{1.73e-77}\\
         \footnotesize{384-bit precision} & \footnotesize{5.08e-116}\\
         \footnotesize{3072-bit precision} & \footnotesize{3.44e-925}\\\hline
    \end{tabular}
    \caption{The machine epsilon of $q$-bit precision floating-point numbers. }
    \label{flo36}
\end{table}

%% file: fjpolynomial.tex
\section{Jacobi fractional polynomials}\label{jfp}

\begin{figure}[!htp]
    \centering
    \subfloat[]{\includegraphics[width=0.31\textwidth]{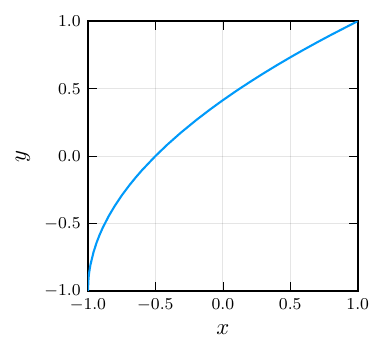}\label{flo20}}
    \quad
    \subfloat[]{\includegraphics[width=0.31\textwidth]{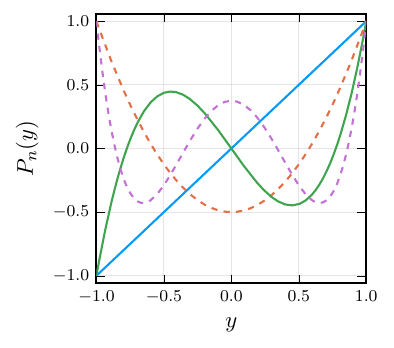}\label{flo21}}
    \quad
    \subfloat[]{\includegraphics[width=0.31\textwidth]{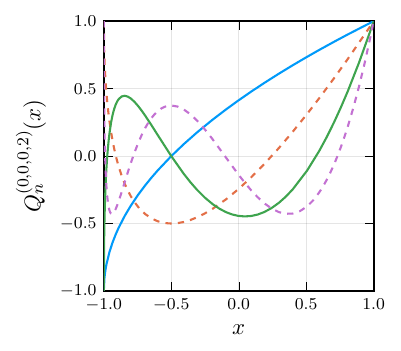}\label{flo22}}
    \caption{ An example of JFP basis functions (\eqref{eq31} with $\alpha = \beta = b =0$ and $p = 2$), which in this case are Legendre polynomials in the mapped variable $y = \sqrt{2(x+1)} - 1$. \protect\subref{flo20}: $\frac{1+x}{2}=\left(\frac{1+y}{2}\right)^2$; \protect\subref{flo21}: $P^{(0,0)}_n(y)$ for $n=1,2,3,4$. \protect\subref{flo22}: $Q^{(0,0,0,2)}_n(x)$ for $n=1,2,3,4$.
    }\label{flo19}
\end{figure}

The JFP basis functions that we shall use are weighted polynomials in fractional powers defined as follows,
\begin{equation}\label{eq31}
    Q_n^{(\alpha,\beta,b,p)}(x):=(1+y)^bP_n^{(\alpha,\beta)}(y), \qquad \frac{1+x}{2}=\left(\frac{1+y}{2}\right)^p, 
\end{equation}
where the choices of $b \in \mathbb{R}$ and $p > 0$ will be problem (and algorithm)-dependent, as we shall find in subsequent sections. The corresponding quasimatrix is denoted by $\mathbf{Q}^{(\alpha,\beta,b,p)}$. Note that the relation between $x$ and $y$ defines a bijection from $[-1, 1]$ to $[-1, 1]$, see \cref{flo20} for an example.  The first few Legendre fractional polynomials\footnote{Recall that Legendre polynomials are the same as Jacobi polynomials with $\alpha = \beta = 0$.} with $b=0, p=2$ are illustrated in \cref{flo21,flo22}.

We use the term Jacobi fractional polynomials, or JFPs, for our basis functions \eqref{eq31} because they define polynomials in fractional powers (or in the mapped variable $y$) and to distinguish them from fractional  Jacobi polynomials~\cite{gogovcheva2005fractional}, which are defined by fractional differential equations. The JFPs are related to the basis functions used by Bhrawy and Zaky by an affine transformation, hence we refer to~\cite{bhrawy2016shifted} for a discussion of some properties of JFPs.

The following result will allow us to relate the JFP basis to a weighted fractional monomial basis in $x$.
\begin{lemma}\label{lem:Myx}
The weighted monomial basis $\mathbf{M}_b(y)$ is equivalent to the following (diagonally scaled) weighted fractional monomial basis in $x$,
\begin{equation}
\mathbf{M}_b(y) = \mathbf{M}^{1/p}_{b/p}(x) D_{b,p}, \label{eq:MyMx}
\end{equation}
where
\begin{align}
D_{b,p}  = 2^{b(1 - 1/p)}\left(\begin{array}{c c c c c}
1 & & & & \\
 & 2^{1 - 1/p} & & & \\
 & & 2^{2(1 - 1/p)} & & \\
  & & & 2^{3(1 - 1/p)} &  \\
   & & & & \ddots
\end{array}    \right)   \label{eq:Dbpdef}.  
\end{align}
\end{lemma}
\begin{proof}
It follows from the definition \eqref{eq:Mgd} and the relation between $x$ and $y$ in \eqref{eq31} that the right-hand side of \eqref{eq:MyMx} simplifies to
\begin{align*}
& \begingroup 
\setlength\arraycolsep{8pt}\left(\begin{array}{c | c | c | c}
2^b\left(\dfrac{1+x}{2}  \right)^{b/p} &  2^{b+1}\left(\dfrac{1+x}{2}  \right)^{(b+1)/p} & 2^{b+2}\left(\dfrac{1+x}{2}  \right)^{(b+2)/p} & \cdots
\end{array}
\right)\endgroup\\
& = \begingroup 
\setlength\arraycolsep{8pt} \left(\begin{array}{c | c | c | c}
(1+y)^b &  (1+y)^{1+b} & (1+y)^{2+b} & \cdots
\end{array}
\right)\endgroup = \mathbf{M}_b(y).
\end{align*}
\end{proof}
\begin{corollary}
The JFP basis is related to a fractional monomial basis as follows
\begin{equation}
\mathbf{Q}^{(\alpha,\beta,b,p)}(x) = \mathbf{M}^{1/p}_{b/p}(x) D_{b,p}  C^{(\alpha,\beta)}. \label{eq:JFPmono}
\end{equation}
\end{corollary}
\begin{proof}
It follows from \eqref{eq6} and  \cref{lem:Myx}  that
\begin{align*}
 \mathbf{Q}^{(\alpha,\beta,b,p)}(x) &=(1+y)^b\mathbf{P}^{(\alpha,\beta)}(y) =(1+y)^b\mathbf{M}_0(y) C^{(\alpha,\beta)} 
             = \mathbf{M}_b(y) C^{(\alpha,\beta)}  \notag \\
            & = \mathbf{M}^{1/p}_{b/p}(x) D_{b,p}  C^{(\alpha,\beta)}. 
\end{align*}
\end{proof}

%% file: fioperator.tex
\section{Properties of fractional integrals applied to the JFP basis}\label{fio}

The results derived in this section will motivate the two algorithms for computing fractional integration matrices that are discussed in the next section. The following result gives the action of $\mathcal{I}^\mu$ on the JFP basis and is the foundation of the first algorithm to be presented in the following section. The subsequent results in this section will inform the second algorithm.

\begin{theorem} \label{prop:fiopex}
Let $\mu, p > 0$ and $\mu p = k$, where $k \in \mathbb{N}_+$, then $\mathcal{I}^\mu  $ maps $\mathbf{Q}^{(\alpha,\beta,b,p)}$ to itself, i.e., $ \mathcal{I}^\mu \mathbf{Q}^{(\alpha,\beta,b,p)} = \mathbf{Q}^{(\alpha,\beta,b,p)} I^{(\alpha,\beta)}_{b,p,\mu}$, where $ I^{(\alpha,\beta)}_{b,p,\mu}$ is a lower banded matrix with $k$ nonzero subdiagonals (i.e., $I^{(\alpha,\beta)}_{b,p,\mu}$ has bandwidths $(k, \infty)$) and it is given by 
\begin{equation}
I^{(\alpha,\beta)}_{b,p,\mu} =  2^{\mu(1 - p)} \left( C^{(\alpha,\beta)}  \right)^{-1}\Lambda^k_{b/p,1/p} C^{(\alpha,\beta)},  \label{eq:fiopc1}
\end{equation}
where $C^{(\alpha,\beta)}$ is defined in \eqref{eq6}, $\left(C^{(\alpha,\beta)}\right)^{-1}$ follows from \eqref{eq:Cabinvdef} and $\Lambda^k_{b/p,1/p}$ is the matrix defined in \eqref{eq:Lkdiag}.
\end{theorem}
\begin{proof}
We derive the result from \eqref{eq:JFPmono} (where $D_{b,p}$ is defined in \eqref{eq:Dbpdef}), \cref{prop:fiomono} with $\delta = b/p$, $\gamma = 1/p$ and the relation between $x$ and $y$ given in \eqref{eq31}:  
\begin{align*}
    \mathcal{I}^\mu \mathbf{Q}^{(\alpha,\beta,b,p)}&=\mathcal{I}^\mu \mathbf{M}^{1/p}_{b/p} D_{b,p} C^{(\alpha,\beta)} 
     = \mathbf{M}^{1/p}_{b/p} \Lambda^k_{b/p,1/p}  D_{b,p} C^{(\alpha,\beta)} \\
    & = \mathbf{Q}^{(\alpha,\beta,b,p)}\left( C^{(\alpha,\beta)}  \right)^{-1}D_{b,p}^{-1}\Lambda^k_{b/p,1/p}D_{b,p} C^{(\alpha,\beta)} \\
     & = 2^{\mu(1 - p)}\mathbf{Q}^{(\alpha,\beta,b,p)}\left( C^{(\alpha,\beta)}  \right)^{-1}\Lambda^k_{b/p,1/p} C^{(\alpha,\beta)},
\end{align*}
where we used the  fact that $D_{b,p}^{-1}\Lambda^k_{b/p,1/p}D_{b,p} = 2^{\mu(1 - p)}\Lambda^k_{b/p,1/p}$ (which follows since the only nonzero diagonal of $\Lambda^k_{b/p,1/p}$ is the $k$-th subdiagonal, see \eqref{eq:Lkdiag})  and hence the result \eqref{eq:fiopc1} follows. The lower banded structure of \eqref{eq:fiopc1} is a consequence of the fact that $ C^{(\alpha,\beta)}$ and  $\left( C^{(\alpha,\beta)}  \right)^{-1}$ are upper triangular and the structure of \eqref{eq:Lkdiag}.
\end{proof}

We know from \cite[Lemma 2.2]{gutleb2020computing} that fractional integrals satisfy recurrence relations in ultraspherical bases. We now derive a general recurrence relation satisfied by fractional integrals, which we shall apply to Jacobi fractional bases.

\begin{proposition}[Fractional integral recurrence]\label{thm1}
    Suppose $u(x)$ is a function on $[-1, 1]$ such that $\mathcal{I}^\mu u(x)$ exists\footnote{For our purposes, it is sufficient to consider functions with convergent expansions of the form $u(x) = (1+x)^{b/p}\sum_{n = 0}^{\infty} c_n(1 + x)^{n/p}$, where $b/p > -1$ and $p \mu  = k \in \mathbb{N}_+$, for which $\mathcal{I}^\mu u(x)$ exists. } for $\mu > 0$, then
    \begin{equation*}
    \mathcal{I}^\mu[\diamond u](x)= x\mathcal{I}^\mu u(x)-\mu\mathcal{I}^{\mu}\mathcal{I}u(x).
    \end{equation*}
\end{proposition}
\begin{proof}
By definition of the fractional integral,
\begin{align*}    
        \mathcal{I}^\mu[\diamond u](x)
            &=\frac{1}{\Gamma(\mu)}\int_{-1}^x\frac{tu(t)}{(x-t)^{1-\mu}}dt\\
            &=\frac{x}{\Gamma(\mu)}\int_{-1}^x\frac{u(t)}{(x-t)^{1-\mu}}dt+\frac{1}{\Gamma(\mu)}\int_{-1}^x\frac{(t-x)u(t)}{(x-t)^{1-\mu}}dt\\
            &=x\mathcal{I}^\mu u(x)-\frac{1}{\Gamma(\mu)}\int_{-1}^x\frac{u(t)}{(x-t)^{-\mu}}dt\\
            &=x\mathcal{I}^\mu u(x)-\mu\mathcal{I}^{\mu+1}u(x)\\
            &=x\mathcal{I}^\mu u(x)-\mu\mathcal{I}^{\mu}\mathcal{I}u(x).
    \end{align*}
\end{proof}

In order to apply \cref{thm1} to the JFP basis, we need to construct matrix representations (which are given in the following lemmata) of multiplication and integration operators. 

\begin{lemma}[Multiplication matrix]\label{thm5}
Suppose $p$ is a positive integer and let $X_{b,p}^{(\alpha,\beta)}$ be the matrix representing multiplication by $x$ in the $\mathbf{Q}^{(\alpha,\beta,b,p)}(x)$ basis, i.e., $x\mathbf{Q}^{(\alpha,\beta,b,p)}(x)=\mathbf{Q}^{(\alpha,\beta,b,p)}(x)X_{b,p}^{(\alpha,\beta)}$, then
\begin{equation*}
X_{b,p}^{(\alpha,\beta)}=2^{1-p}\left(L_{(\alpha,\beta+1)}^{(\alpha,\beta)}R_{(\alpha,\beta)}^{(\alpha,\beta+1)}\right)^p-\mathit{1},
\end{equation*}
 where $\mathit{1}$ denotes the identity matrix, and $X_{b,p}^{(\alpha,\beta)}$ has bandwidths $(p, p)$.
\end{lemma}
\begin{proof}
The result follows from \eqref{eq31} and \cref{thm2}:
\begin{align*}
x\mathbf{Q}^{(\alpha,\beta,b,p)}(x) &= (1+y)^b\left[ 2^{1-p}\left(1 + y  \right)^p - 1  \right]\mathbf{P}^{(\alpha,\beta)}(y) \\
& = (1+y)^b\mathbf{P}^{(\alpha,\beta)}(y) \left(2^{1-p} \left(L_{(\alpha,\beta+1)}^{(\alpha,\beta)}R_{(\alpha,\beta)}^{(\alpha,\beta+1)}\right)^p-\mathit{1}  \right) \\
& = \mathbf{Q}^{(\alpha,\beta,b,p)}(x) \left( 2^{1-p}\left(L_{(\alpha,\beta+1)}^{(\alpha,\beta)}R_{(\alpha,\beta)}^{(\alpha,\beta+1)}\right)^p-\mathit{1}  \right).
\end{align*}
Since $L_{(\alpha,\beta+1)}^{(\alpha,\beta)}R_{(\alpha,\beta)}^{(\alpha,\beta+1)}$ has bandwidths $(1, 1)$, it follows that  $X_{b,p}^{(\alpha,\beta)}$ has bandwidths $(p, p)$.
\end{proof}

\begin{lemma}[Integration matrix]\label{thm4}
 Suppose $\{p,\beta-b,b+p-1-\beta\}\subset\mathbb{N}_0$,  and denote the matrix representing the integral operator applied to  the JFP basis by $I_{b,p}^{(\alpha,\beta)}$, i.e.,  $\mathcal{I}[\mathbf{Q}^{(\alpha,\beta,b,p)}]=\mathbf{Q}^{(\alpha,\beta,b,p)}I_{b,p}^{(\alpha,\beta)}$, then
\begin{equation}\label{eq:iop}
    I_{b,p}^{(\alpha,\beta)}=\frac{p}{2^{p-1}}L_{(\alpha,\beta+p)}^{(\alpha,\beta)}R_{(\alpha-1,b+p)}^{(\alpha,\beta+p)}\left(W^{(\alpha,b+p-1)}_{(\alpha-1,b+p)}\right)^{-1}R_{(\alpha,\beta)}^{(\alpha,b+p-1)},
    \end{equation}
    and $I_{b,p}^{(\alpha,\beta)}$ has bandwidths $(p, p)$.
\end{lemma}
\begin{proof}
Using \eqref{eq:wintop}, the matrices in \cref{flo1}, the definition of the JFP basis in \eqref{eq31} and making the change of variables $\frac{1 + t}{2} = \left(\frac{1+s}{2}\right)^p$, it follows that
\begin{align*}    
     \mathcal{I}[\mathbf{Q}^{(\alpha,\beta,b,p)}](x) & =   \int_{-1}^x\mathbf{Q}^{(\alpha,\beta,b,p)}(t)dt \\
            &=\frac{p}{2^{p-1}}\int_{-1}^y(1+s)^{b+p-1}\mathbf{P}^{(\alpha,\beta)}(s)ds\\
            &=\frac{p}{2^{p-1}}\mathcal{I}[(1+\diamond)^{b+p-1}\mathbf{P}^{(\alpha,\beta)}](y)\\
            &=\frac{p}{2^{p-1}}\mathcal{I}[(1+\diamond)^{b+p-1}\mathbf{P}^{(\alpha,b+p-1)}](y)R_{(\alpha,\beta)}^{(\alpha,b+p-1)}\\
            &=\frac{p}{2^{p-1}}(1+y)^{b+p}\mathbf{P}^{(\alpha-1,b+p)}(y)\left(W^{(\alpha,b+p-1)}_{(\alpha-1,b+p)}\right)^{-1}R_{(\alpha,\beta)}^{(\alpha,b+p-1)},
\end{align*}
where
\begin{align*}
	(1+y)^{b+p}\mathbf{P}^{(\alpha-1,b+p)}(y) & = (1+y)^{b+p}\mathbf{P}^{(\alpha,\beta+p)}(y)R_{(\alpha-1,b+p)}^{(\alpha,\beta+p)} \\
	&=(1+y)^{b}\mathbf{P}^{(\alpha,\beta)}(y)L_{(\alpha,\beta+p)}^{(\alpha,\beta)}R_{(\alpha-1,b+p)}^{(\alpha,\beta+p)}\\
	&=\mathbf{Q}^{(\alpha,\beta,b,p)}(x)L_{(\alpha,\beta+p)}^{(\alpha,\beta)}R_{(\alpha-1,b+p)}^{(\alpha,\beta+p)},
\end{align*}
hence we reach \eqref{eq:iop}.

According to \cref{flo1}, the bandwidths of $L_{(\alpha,\beta+p)}^{(\alpha,\beta)}$, $R_{(\alpha-1,b+p)}^{(\alpha,\beta+p)}$ and $R_{(\alpha,\beta)}^{(\alpha,b+p-1)}$ are $(p,0)$, $(0,\beta-b+1)$ and $(0, b+p-1-\beta)$, respectively. Recalling that $W^{(\alpha,b+p-1)}_{(\alpha-1,b+p)}$ is diagonal, we conclude that $I_{b,p}^{(\alpha,\beta)}$ has bandwidths $(p,p)$.
\end{proof}

\begin{theorem}\label{thm:sylveqs}
Suppose $\{p,\beta-b,b+p-1-\beta\}\subset\mathbb{N}_0$, then the fractional integration matrix defined in \cref{prop:fiopex} satisfies the following banded Sylvester equations
\begin{align}
I_{b,p,\mu}^{(\alpha,\beta)}I_{b,p}^{(\alpha,\beta)} & =I_{b,p}^{(\alpha,\beta)}I_{b,p,\mu}^{(\alpha,\beta)}, \label{eq:intcomm} \\
    I_{b,p,\mu}^{(\alpha,\beta)}\left(X_{b,p}^{(\alpha,\beta)}+\mu I_{b,p}^{(\alpha,\beta)}\right) & =X_{b,p}^{(\alpha,\beta)}I_{b,p,\mu}^{(\alpha,\beta)}, \label{eq18} 
\end{align}
where $X_{b,p}^{(\alpha,\beta)}$ and $I_{b,p}^{(\alpha,\beta)}$  have bandwidths $(p,p)$.
\end{theorem}
\begin{proof}
The bandwidths of the matrices follow from \cref{thm5,thm4}. The first equation \eqref{eq:intcomm} follows immediately from the commutativity of fractional (and integer-order) integration matrices stated in \eqref{eq20}. To derive \eqref{eq18}, we apply \cref{thm1} to the JFP basis, then
\begin{align}
\mathcal{I}^\mu\left[ \diamond \mathbf{Q}^{(\alpha,\beta,b,p)}\right](x)  & = \mathcal{I}^\mu \mathbf{Q}^{(\alpha,\beta,b,p)}(x)X_{b,p}^{(\alpha,\beta)} \notag \\
& =\mathbf{Q}^{(\alpha,\beta,b,p)}(x)I_{b,p,\mu}^{(\alpha,\beta)}X_{b,p}^{(\alpha,\beta)} \label{eq:Iopeqlhs} \\
& = x \mathcal{I}^\mu \mathbf{Q}^{(\alpha,\beta,b,p)}(x) - \mu \mathcal{I}^\mu \mathcal{I} \mathbf{Q}^{(\alpha,\beta,b,p)}(x) \notag \\
& = x \mathbf{Q}^{(\alpha,\beta,b,p)}(x) I_{b,p,\mu}^{(\alpha,\beta)} - \mu \mathcal{I}^\mu \mathbf{Q}^{(\alpha,\beta,b,p)}(x) I_{b,p}^{(\alpha,\beta)} \notag \\
& =    \mathbf{Q}^{(\alpha,\beta,b,p)}(x) X_{b,p}^{(\alpha,\beta)}  I_{b,p,\mu}^{(\alpha,\beta)} - \mu \mathbf{Q}^{(\alpha,\beta,b,p)}(x)  I_{b,p,\mu}^{(\alpha,\beta)} I_{b,p}^{(\alpha,\beta)} \label{eq:Iopeqrhs},
\end{align}
comparing \eqref{eq:Iopeqrhs} and \eqref{eq:Iopeqlhs}, \eqref{eq18} follows.
\end{proof}

\section{Algorithms for computing fractional integration matrices}\label{sect:fioalgs}
 In this section we discuss the following algorithms for computing the entries of fractional integral matrices:
\begin{itemize}
\item[] \textbf{Algorithm 1:} Compute $I^{(\alpha,\beta)}_{b,p,\mu}$ by solving (see \cref{prop:fiopex})
\begin{equation}
C^{(\alpha,\beta)}  I^{(\alpha,\beta)}_{b,p,\mu} =  2^{\mu(1 - p)} \Lambda^k_{b/p,1/p} C^{(\alpha,\beta)}. \label{eq:fiosyst}
\end{equation}
\item[] \textbf{Algorithm 2:} Compute $I^{(\alpha,\beta)}_{b,p,\mu}$ by solving either of the banded Sylvester equations \eqref{eq:intcomm} or \eqref{eq18}. This requires the pre-computation of the first $p$ columns of $I^{(\alpha,\beta)}_{b,p,\mu}$, which can be obtained from Algorithm 1.
\end{itemize}

\subsection{Algorithm 1}

In Algorithm 1, each column of $I^{(\alpha,\beta)}_{b,p,\mu}$ is computed independently (hence making Algorithm 1 amenable to parallelisation) by solving an upper triangular linear system (recall that $C^{(\alpha,\beta)}$ is upper triangular). If $I^{(\alpha,\beta)}_{b,p,\mu}$ has bandwidths $(k,\infty)$, $k \geq 0$, 
then obtaining the $j$-th column (with $j \geq 0$) requires the solution of a $(j+1+k)\times(j+1+k)$ system. Hence, computing $N$ columns of $I^{(\alpha,\beta)}_{b,p,\mu}$  has a complexity of $\mathcal{O}(N^3)$. However, the condition numbers of $N \times N$ sections of $C^{(\alpha,\beta)}$ (which we denote by $C^{(\alpha,\beta)}_{0:N-1,0:N-1}$) grow exponentially~\cite{gautschi1979condition} with $N$ and thus high-precision arithmetic is required to maintain accuracy. In   $q$-bit precision, each multiplication costs $\mathcal{O}(q\log q\log\log q)$ operations with the Schönhage--Strassen algorithm\footnote{There are galactic algorithms with lower asymptotic complexity \cite{furer2009faster,harvey2021integer}.},  hence Algorithm~1 has a complexity of $\mathcal{O}(qN^3\log q \log \log q)$.  As shown in \cref{ps}, $q = \mathcal{O}(N)$ is required to maintain accuracy due to the exponentially growing condition numbers of $C^{(\alpha,\beta)}$, there a pseudo-stabilized version of Algorithm~1 has a complexity of $\mathcal{O}(N^4\log N \log\log N)$.    For Algorithm 1, we may choose any $p>0$ (rational or irrational), such that $\mu p = k \in \mathbb{N}_+$. 
 
\begin{figure}[!htp]
    \centering
    \subfloat[Solving \eqref{eq:fiosyst} via backslash]{\includegraphics[width=0.45\textwidth]{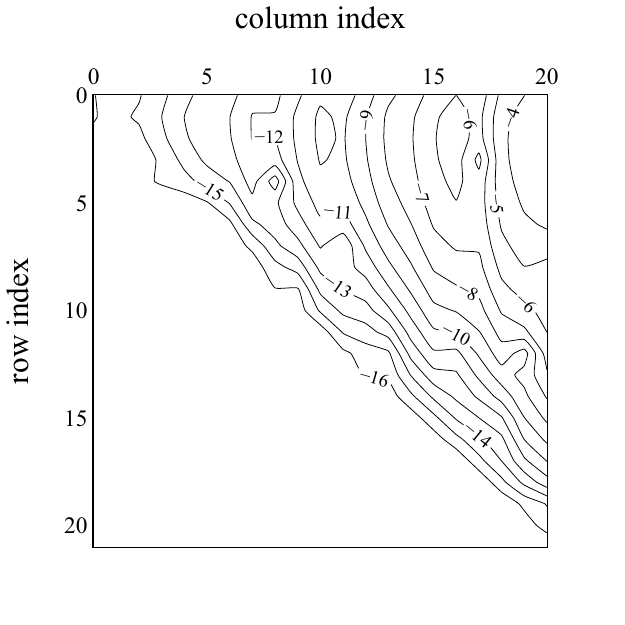}}
    \quad
    \subfloat[Solving \eqref{eq:fiosyst} via the explicit inverse]{\includegraphics[width=0.45\textwidth]{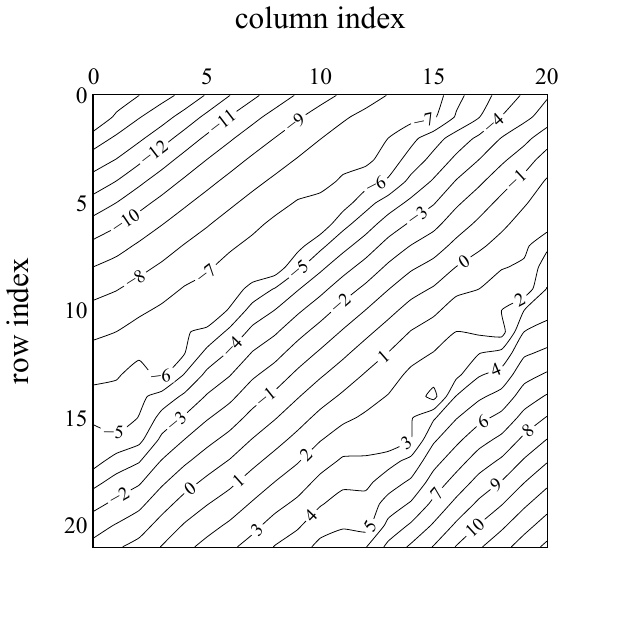}}
    \caption{Errors in $\log_{10}$ scale from using Algorithm 1 to compute $I_{0,2,1/2}^{(0,0)}$.}\label{flo8}
\end{figure}

In \cite{bhrawy2016shifted}, \eqref{eq:fiosyst} was solved by forming the inverse of $C^{(\alpha,\beta)}$ explicitly, using \eqref{eq:Cabinvdef}. However, this is more unstable (see \cref{flo8}) and not faster than solving \eqref{eq:fiosyst} via the backslash ($\backslash$) command in \texttt{Julia}.

\subsection{Algorithm 2}\label{sect:alg2sect}

The Sylvester equations  \eqref{eq:intcomm} and \eqref{eq18} take the form
\begin{equation*}
AB = CA,
\end{equation*}
where $A = I^{(\alpha,\beta)}_{b,p,\mu}$ and $B$ and $C$ are banded matrices with bandwidths $(p, p)$. Considering the $(m, n)$-entry of both sides of the equation, and using the bandedness of $B$ and $C$, it follows that
\begin{equation*}
A_{m,n-p:n+p}B_{n-p:n+p,n}=C_{m,m-p:m+p}A_{m-p:m+p,n},
\end{equation*}
where an entry is considered to be zero if an index is negative, and thus
\begin{equation}\label{eq19}
   A_{m,n+p}B_{n+p,n}=C_{m,m-p:m+p}A_{m-p:m+p,n}-A_{m,n-p:n+p-1}B_{n-p:n+p-1,n}. 
\end{equation}
This shows that if the first $p$ columns of $A$ are known (which, for $A = I^{(\alpha,\beta)}_{b,p,\mu}$, can be obtained from Algorithm 1), then the subsequent columns can be obtained from the recurrence relation defined by \eqref{eq19}.

\begin{figure}[!htp]
    \centering
    \subfloat[Solution to \eqref{eq18}]{\includegraphics[width=0.45\textwidth]{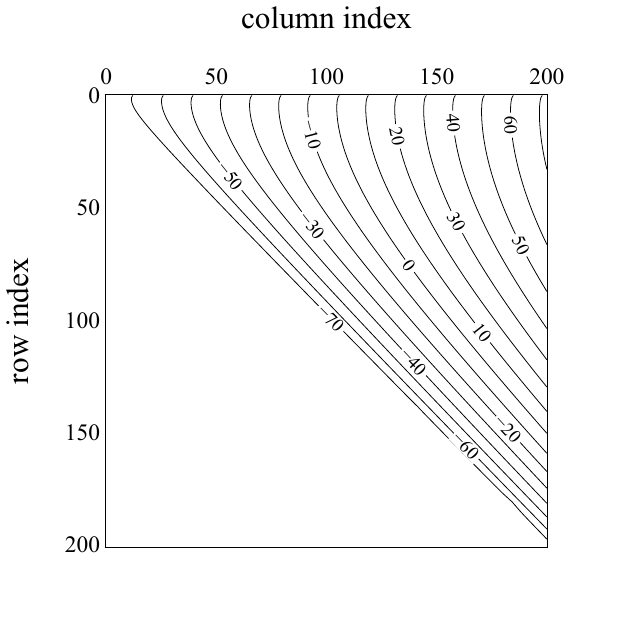}}
    \quad
    \subfloat[Solution to \eqref{eq:intcomm}]{\includegraphics[width=0.45\textwidth]{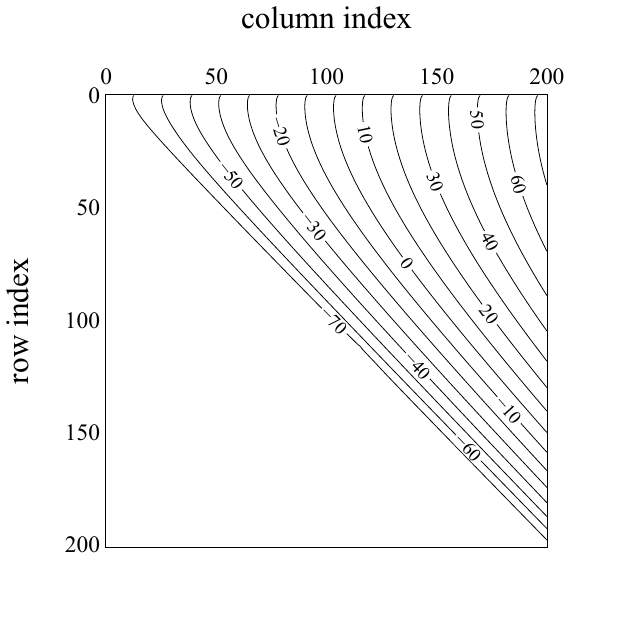}}
    \caption{Errors in $\log_{10}$ scale from computing $I_{0,2,1/2}^{(0,0)}$ by solving the banded Sylvester equations \eqref{eq18} and \eqref{eq:intcomm} in 256-bit precision. The difference in accuracy between the two solutions is negligible.}\label{flo9}
\end{figure}

We have found that Algorithms~1 and 2 are roughly equally unstable. In Algorithm~2, the numerical errors propagate due to column-by-column recurrence, and thus the accuracy is limited by the precision we begin with. By contrast, in Algorithm~1, each column can be computed independently to any desired precision. Another disadvantage of Algorithm~2 is that we can only use it for rational $\mu$ since we require $\mu p = k \in \mathbb{N}_+$ and $p$ must be a positive integer (see \cref{thm:sylveqs}).
However, the complexity of Algorithm~2 is lower than that of Algorithm~1 by one order. Specifically, the first $N$ columns of $I^{(\alpha,\beta)}_{b,p,\mu}$ involve the  computation of $\mathcal{O}(N^2)$ entries, each of which requiring up to $4p+3$ multiplications, leading to a total cost of $\mathcal{O}(pN^2q\log q\log\log q)$ in $q$-bit precision. As shown in~\cref{ps}, $q=\mathcal{O}(N)$ and we conclude that Algorithm~2 has an overall complexity of $\mathcal{O}(pN^3\log N\log\log N)$.
\cref{flo9} illustrates the growth of errors for a high-precision computation of $I^{(\alpha,\beta)}_{b,p,\mu}$ obtained via Algorithm~2.

In \cref{ps} we describe how we pseudostabilize Algorithm~2 by automatically incorporating high-precision arithmetic in such a manner that the first $N$ columns of the matrix $I^{(\alpha,\beta)}_{b,p,\mu}$ are  computed to a prescribed accuracy while minimizing computational cost. 

Algorithm~2 can be stabilised with Algorithm~1 at the expense of higher computational complexity as follows. We can use Algorithm~1 to compute the first $p$ columns of $I^{(\alpha,\beta)}_{b,p,\mu}$ and also its $N$-th column for some $N > p$ with high-precision arithmetic, which costs $\mathcal{O}(N^3\log N\log\log N )$ operations. Then the entries in the intervening columns (columns $p+1$ to $N-1$) can be assembled into a single vector and approximated as the least squares solution to an overdetermined system. We have found that this system is well-conditioned and can thus be computed accurately in a lower precision (for example, we have achieved roughly $10^{-13}$ accuracy in double precision). However, we found that this procedure is slower than the pseudo-stabilized Algorithm~2. 

Another strategy to stabilize Algorithm~2, which we shall pursue in future work, is to derive asymptotic approximations to the entries of $I^{(\alpha,\beta)}_{b,p,\mu}$ in column $N$ as $N \to \infty$, which could be used (as described above) to set up an overdetermined system and stably compute the entries of  $I^{(\alpha,\beta)}_{b,p,\mu}$. This approach would have the optimal  $\mathcal{O}(N^2)$ complexity.

%% file: spectral.tex
\section{Spectral methods using Jacobi fractional polynomials}\label{spec}

We first test our methods  by applying them to problems that have solutions expressible in terms of Mittag--Leffler functions in \cref{sect:testprobs} before tackling more challenging problems in \cref{sect:probs}. In \cref{exa3} of \cref{sect:testprobs} we shall consider a striking example in which the performance of the JFP method and the sum space method of~\cite{hale2018fast} are radically different.

\subsection{FIEs, FDEs and a fractional PDE with Mittag--Leffler function solutions}\label{sect:testprobs}

The general Mittag--Leffler functions are defined as
\begin{equation}\label{eq24}
    E_{\alpha,\beta}(z)=\sum_{k=0}^\infty\frac{z^k}{\Gamma(\beta+\alpha k)},\qquad \alpha,\beta,z\in\mathbb{C},\qquad  \mathrm{Re\,}\alpha,\mathrm{Re\,}\beta>0,
\end{equation}
and they are the general solutions to the following fractional integral equations (which can be verified using \eqref{eq4}):
\begin{theorem}[{\cite[Theorem 13.4]{haubold2011mittag}}]\label{thm6}
    If $\mathrm{Re\,}\nu,\mathrm{Re\,}\mu >0$, then the solution to the fractional integral equation
    \begin{equation}
        u(x)+\lambda\mathcal{I}^\mu u(x)= (1+x)^{\nu-1}, \qquad x \in [-1, 1],  \label{eq:MLeq}
    \end{equation}
    is given by\footnote{In \cite[Theorem 13.4]{haubold2011mittag} there is a factor of $\Gamma(\nu)$ in the formula for $u(x)$, which we believe is due  to a miscalculation.}
    \begin{equation}
        u(x)=(1+x)^{\nu-1}E_{\mu,\nu}(-\lambda(1+x)^\mu).  \label{eq:MLsol}
    \end{equation}
\end{theorem}

For the problems we consider, $\alpha$ and $\beta$ are positive real numbers, in which case $E_{\alpha,\beta}(z)$ is an entire function. Various methods for computing Mittag--Leffler functions have been implemented, see \cite{garrappa2015numerical,gorenflo2002computation,hilfer2006computation,ortigueira2019numerical}. We shall use~\cite{garrappa2015numerical} to benchmark the accuracy of our methods.

We can express the solution \eqref{eq:MLsol} in the normalised monomial basis \eqref{eq:normM0} as $u(x) = (1+x)^{\nu-1}\widetilde{\mathbf{M}}_{0}^{\mu}(x) \widetilde{\bm{u}}_{\nu,\mu} $, where 
\begin{equation}
\widetilde{\bm{u}}_{\nu,\mu} = \begingroup 
\setlength\arraycolsep{8pt} \left(\begin{array}{c c c c c}
\displaystyle{\frac{1}{\Gamma(\nu)}} & \displaystyle{\frac{-2^{\mu}\lambda }{\Gamma(\nu + \mu)}} & \displaystyle{\frac{(-2^{\mu}\lambda )^2 }{\Gamma(\nu+2\mu)}} & \displaystyle{\frac{(-2^{\mu}\lambda )^3 }{\Gamma(\nu+3\mu)}} & \cdots 
\end{array} \right)^{\top}.\endgroup 
 \label{eq:monocoeffs}
\end{equation}

If we define $\bm{u}_{\alpha,\beta,b,p}$ to be the coefficients of the solution in the JFP basis, i.e., $u(x) = \mathbf{Q}^{(\alpha,\beta,b,p)}(x) \bm{u}_{\alpha,\beta,b,p}$, then using \eqref{eq:JFPmono}, it follows that
\begin{equation}
u(x) = (1+x)^{\nu-1}\widetilde{\mathbf{M}}_{0}^{\mu}(x) \widetilde{\bm{u}}_{\nu,\mu} = \mathbf{Q}^{(\alpha,\beta,b,p)}(x) \bm{u}_{\alpha,\beta,b,p} = 
\mathbf{M}^{1/p}_{b/p}(x)\bm{u}_{b,p},  \label{eq:urep}
\end{equation}
where
\begin{equation}
D_{b,p}C^{(\alpha,\beta)}\bm{u}_{\alpha,\beta,b,p} = \bm{u}_{b,p}.  \label{eq:MLcoeffs}
\end{equation}
For $\mathcal{I}^\mu u$  to be defined in the JFP basis, a necessary condition is $b > -p$, to ensure that the basis  $\mathbf{Q}^{(\alpha,\beta,b,p)}$ is integrable (because it contains monomials of the form $(1+x)^{b/p}$).
Equations \eqref{eq:urep} and \eqref{eq:MLeq} therefore hold only if there are integers $k_* \geq 1$ and $n \geq 0$ such that
\begin{equation}
p  = \frac{k_*}{\mu}>-b,\qquad  \frac{b+n}{p} = \nu - 1. \label{eq:pbvals}
\end{equation}
If we let $\bm{u}(\ell)$ denote the $\ell$-th entry of the vector $\bm{u}$ for $\ell \geq 0$, then $\bm{u}_{b,p}(\ell)$ is the coefficient of the monomial basis function $(1+x)^{(b+\ell)/p}$. Likewise, $\widetilde{\bm{u}}_{\nu,\mu}(\ell)$ is the coefficient of $(1+x)^{\nu-1+\ell\mu}$.  Hence, \eqref{eq:urep} and \eqref{eq:pbvals} imply that $\bm{u}_{b,p}(n +mk_*) =\widetilde{\bm{u}}_{\nu,\mu}(m)$, $m \geq 0$, and otherwise, $\bm{u}_{b,p}(\ell) = 0$, $\ell \neq n +mk_*$. Specifying $\bm{u}_{b,p}$ thus, we can obtain $\bm{u}_{\alpha,\beta,b,p}$ by solving \eqref{eq:MLcoeffs}, which is an ill-conditioned system whose condition numbers grow exponentially and relies on knowing the exact coefficients of the solution in the monomial basis. Instead, we shall obtain the JFP coefficients $\bm{u}_{\alpha,\beta,b,p}$ by solving the integral equation \eqref{eq:MLeq} in the JFP basis which, as we shall see, result in systems whose condition numbers are bounded by a quantity that grows linearly in $\lambda$.

If the conditions \eqref{eq:pbvals} are satisfied and we set $u(x) = \mathbf{Q}^{(\alpha,\beta,b,p)}(x) \bm{u}_{\alpha,\beta,b,p}$, then \eqref{eq:MLeq} becomes
\begin{equation*}
\mathbf{Q}^{(\alpha,\beta,b,p)}\left( \mathit{1} + \lambda I^{(\alpha,\beta)}_{b,p,\mu}  \right)     \bm{u}_{\alpha,\beta,b,p} =  \mathbf{Q}^{(\alpha,\beta,b,p)}\bm{f},
\end{equation*}
where $(1 + x)^{\nu - 1} = \mathbf{Q}^{(\alpha,\beta,b,p)}(x)\bm{f}$, i.e., $\bm{f}$ contains the coefficients of $(1 + x)^{\nu - 1}$ in the JFP basis. Using  \eqref{eq31} and the fact that $(b+n)/p = \nu - 1$, it follows that $(1 + x)^{\nu - 1} = 2^{(b+n)(1/p-1)}(1+y)^{b+n}$, hence $\bm{f}$ has $n+1$ nonzero entries.   In the simplest case, $n = 0$ and $\bm{f} = f_0\bm{e}_0$, where $\bm{f}(0) = f_0$. By \cref{prop:fiopex}, the matrix $I^{(\alpha,\beta)}_{b,p,\mu}$ has bandwidths $(k_*,\infty)$ and we obtain $\bm{u}_{\alpha,\beta,b,p}$ by truncating and solving the system
\begin{equation}
\left( \mathit{1} + \lambda I^{(\alpha,\beta)}_{b,p,\mu}  \right)     \bm{u}_{\alpha,\beta,b,p} =  \bm{f}.  \label{eq:MLeqsyst}
\end{equation}

If \eqref{eq:pbvals} is satisfied and we use the relation between $x$ and $y$ in \eqref{eq31}, then  the exact solution \eqref{eq:MLsol} can be expressed as
\begin{equation}
u(x) = 2^{(b+n)(1/p-1)}(1 + y)^{b+n}E_{\mu,\nu}\left(-\lambda 2^{k_*(1/p-1)}(1 + y)^{k_*} \right),  \label{eq:exmappedsol}
\end{equation}
which we represent in the JFP basis as 
\begin{equation*}
u(x) = \mathbf{Q}^{(\alpha,\beta,b,p)}(x)\bm{u}_{\alpha,\beta,b,p} = (1 + y)^b \mathbf{P}^{(\alpha,\beta)}(y)\bm{u}_{\alpha,\beta,b,p}.
\end{equation*}
Comparing these two equations, we conclude that by truncating and solving \eqref{eq:MLeqsyst}, we are approximating an entire function of $y$ by a truncated expansion in Jacobi polynomials. Therefore, we can use bounds on Jacobi coefficients of Jacobi expansions of analytic functions, which are given in~\cite{zhao2013sharp}, to obtain bounds\footnote{For now, we do not take into account the error incurred by solving the system \eqref{eq:MLeqsyst}. This will require an analysis of the conditioning of the system, which we leave for future work. Numerically, we shall find in \cref{FLO18} that the condition number of $N \times N$  truncations of \eqref{eq:MLeqsyst} are bounded as $N \to \infty$ by a constant that grows as $\mathcal{O}(\lambda)$.} on the entries of $\bm{u}_{\alpha,\beta,b,p}$. For example, if we use Chebyshev polynomials, then
\begin{equation}
\vert \bm{u}_{\alpha,\beta,b,p}(n) \vert \leq 2M \rho^{-n}, \qquad n \geq 0, \qquad  \alpha = \beta = -1/2,  \label{eq:chebcoeffsb}
\end{equation}   
holds for any $\rho > 1$ and thus the coefficients decay super-exponentially~\cite[Chapter 8]{trefethen2019approximation}. Bounds similar to \eqref{eq:chebcoeffsb} hold for any Jacobi parameters $\alpha, \beta > -1$, see~\cite{zhao2013sharp}.

\begin{remark}
In \eqref{eq:chebcoeffsb}, $\rho$  refers to the sum of the semi-axes of a Bernstein ellipse $E_{\rho}$ (in the complex $y$-plane) on and within which the function we are approximating (the Mittag--Leffler function) is analytic and $M$ is the maximum modulus of the function on $E_{\rho}$.
\end{remark}

We shall find it instructive to compare the bounds on the JFP coefficients in \eqref{eq:chebcoeffsb} to an estimate of the largest coefficient in the (fractional) normalized monomial basis, i.e., $\max_{n \geq 0} \vert \widetilde{\bm{u}}_{\nu,\mu}(n) \vert$, where (see \eqref{eq:monocoeffs}),
\begin{equation}
\widetilde{\bm{u}}_{\nu,\mu}(n) = \frac{(-2^{\mu}\lambda )^n}{\Gamma(\nu + n \mu)}, \qquad n \geq 0.  \label{eq:mutilde}
\end{equation}
The large-$\lambda$ regime is of particular interest since it will arise in the solution to a fractional heat/wave equation. 

We let $0 < \mu < 1$ and consider the case $\lambda \gg 1$. As $n$ increases from $0$, the magnitude of $\widetilde{\bm{u}}_{\nu,\mu}(n)$ increases and reaches a maximum at which
\begin{equation}
 \vert \widetilde{\bm{u}}_{\nu,\mu}(n+1) \vert \approx \vert \widetilde{\bm{u}}_{\nu,\mu}(n) \vert, \qquad \Rightarrow \qquad \frac{\Gamma(\nu + (n+1)\mu)}{\Gamma(\nu + n \mu)} \approx 2^\mu \lambda, \label{eq:maxcond}
\end{equation} 
then, as $n$ increases further, $\vert \widetilde{\bm{u}}_{\nu,\mu}(n) \vert$ decays super-exponentially as $n \to \infty$. If $\lambda$ is large, then \eqref{eq:maxcond} is satisfied for large $n$ and we can use a leading order approximation of the ratio of gamma functions in \eqref{eq:maxcond} from~\cite[Eq.~(1)]{tricomi1951asymptotic} as $n \to \infty$,
\begin{equation*}
\frac{\Gamma(\nu + (n+1)\mu)}{\Gamma(\nu + n \mu)} \sim  (n\mu)^\mu \approx 2^\mu \lambda, \qquad \lambda \to \infty.
\end{equation*}   
Using the approximation $(n\mu)^\mu \approx 2^\mu \lambda$, $\Gamma(\nu + n \mu ) \sim (n\mu)^\nu\Gamma(n \mu )$, $n \to \infty$ (which again follows from the result in~\cite{tricomi1951asymptotic}) and Stirling's approximation in \eqref{eq:mutilde}, it follows that
\begin{equation}
\max_{n \geq 0} \vert \widetilde{\bm{u}}_{\nu,\mu}(n) \vert \approx \frac{\left( 2\lambda^{1/\mu} \right)^{1/2-\nu}}{\sqrt{2\pi}}\exp\left(2\lambda^{1/\mu}\right), \qquad \lambda \gg 1.  \label{eq:maxmonoc}
\end{equation}
Since $0 < \mu < 1$, this shows that the largest monomial coefficient grows super-exponentially in $\lambda \gg 1$. We shall find evidence of super-exponential growth (as a function of $\lambda$) in the maximum coefficient of the sum space solution in \cref{exa3}, which renders the method catastrophically unstable for large $\lambda$.  By contrast, we shall find in this same example that the JFP coefficients (with bounds given in \eqref{eq:chebcoeffsb}) are bounded below $1$ for all the values of $\lambda$ that we consider.

\begin{example}\label{exa1}
We first consider a simple instance of \eqref{eq:MLeq}, viz.,
\begin{equation}
u(x)+\mathcal{I}^\mu u(x)=1, \qquad x \in [-1, 1],  \label{eq:exFIE1}
\end{equation} 
which has the exact solution $u(x)=E_{\mu,1}\left(-(1+x)^\mu\right)$, according to \cref{thm6}. The conditions \eqref{eq:pbvals} with $\nu = 1$ imply that we require $b = -n$ and $b > -p$. Hence, $b$ can take any one of the following integer values: $b = -\lceil p-1 \rceil, \ldots, -1, 0$. If $\mu$ is rational, we can let $p$ be an integer and choose $\beta$ such that the conditions $\beta - b, p-1-(\beta-b) \in \mathbb{N}_0$ in \cref{thm:sylveqs} are satisfied, which enables us to use Algorithm 2 (see \cref{sect:fioalgs}).

\cref{FLO3} illustrates that truncating and solving \eqref{eq:MLeqsyst} converges super-exponentially for a wide range of fractional orders (including irrational orders), as expected. We let $\alpha=\beta=b=0$, vary $\mu$ and $p$ and fix the value of  $k_* = 1$ (i.e., $\mu p = 1$) in \cref{FLO24,FLO25}  and in \cref{FLO26,FLO27}, we let $k_* \geq 1$. \cref{FLO26,FLO27} suggest that the choice $k_* = 1$ yields faster convergence than choosing $k_* > 1$. However, for some equations (e.g.,  \eqref{eq:exFIE3}), larger $k_*$ can have (sometimes significantly) better performance. We shall pursue the question of determining optimal choices of parameter values (which might involve optimising the bounds \eqref{eq:chebcoeffsb}) in future work.

\begin{figure}[!ht]
    \centering
    \subfloat[$\mu=1/p$ for different $p$]{\includegraphics[width=0.4\textwidth]{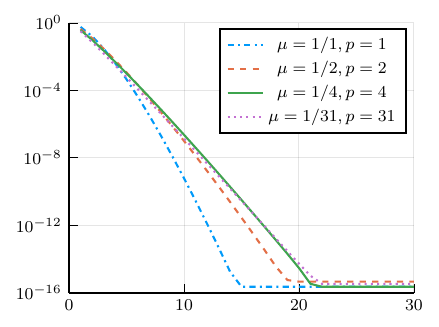}\label{FLO24}}
    \quad
    \subfloat[$p=1/\mu$ for different irrational $\mu$]{\includegraphics[width=0.4\textwidth]{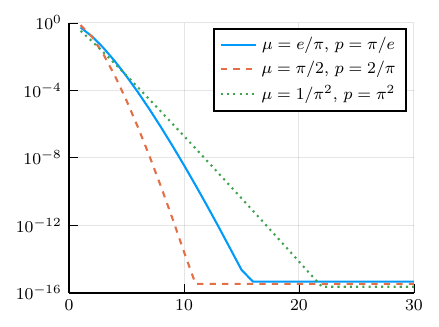}\label{FLO25}}
    
    \subfloat[Different $\mu$ for $p=3$]{\includegraphics[width=0.4\textwidth]{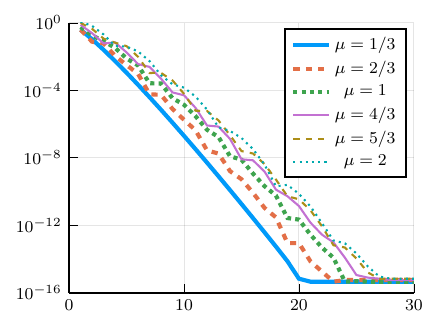}\label{FLO26}}
    \quad
    \subfloat[Different $p$ for $\mu=1/2$]{\includegraphics[width=0.4\textwidth]{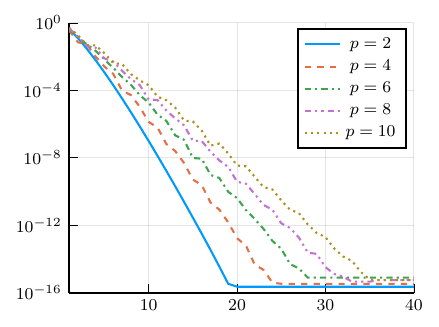}\label{FLO27}}
    \caption{Errors produced by the JFP method for the problem \eqref{eq:exFIE1} for different fractional orders with $\alpha=\beta=b=0$. The vertical axes are the maximum errors (calculated by evaluating the computed solution on the grid $[-1\!\!:\!\!0.01\!\!:\!\!1]$) and the horizontal axes give the truncation size of the system \eqref{eq:MLeqsyst} that is solved to compute the JFP coefficients.}\label{FLO3}
\end{figure}
\end{example}

\begin{example}\label{exa2}
We modify the previous example by letting the right-hand side have a singularity,
\begin{equation}
u(x)+\mathcal{I}^{1/3} u(x)=\sqrt{1+x}, \qquad x \in [-1,1].\label{eq:exFIE2}
\end{equation}
By \cref{thm6}, the exact solution is $u(x)=\sqrt{1+x}E_{\frac{1}{3},\frac{3}{2}}\left(-(1+x)^{1/3}\right)$, which is shown in \cref{FLO39}. We let $k_* = 1$ in \eqref{eq:pbvals}, then $p = 3$ and $b$ is constrained to the values $b = 3/2 - n$, with $b >-p=-3$ and $n \geq 0$, hence the permissible values are $b = 3/2, 1/2, -1/2, -3/2, -5/2$. To use Algorithm~2, we also require $\beta-b\in\{0,1,2\}$, see \cref{thm:sylveqs}. We let $\alpha=\beta=b=-1/2$ and the resulting convergence of the numerical solutions is shown in \cref{FLO37}.

\begin{figure}[!htp]
    \centering  
    \subfloat[The exact solution]{\includegraphics[width=0.4\textwidth]{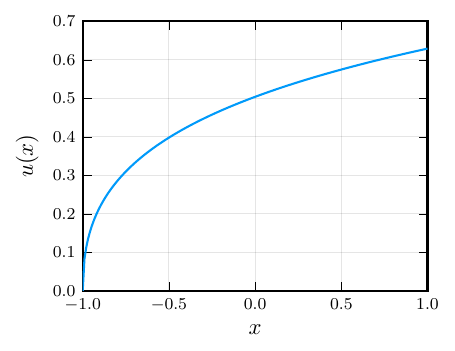}\label{FLO39}}
    \subfloat[Error and decay of coefficients]{\includegraphics[width=0.4\textwidth]{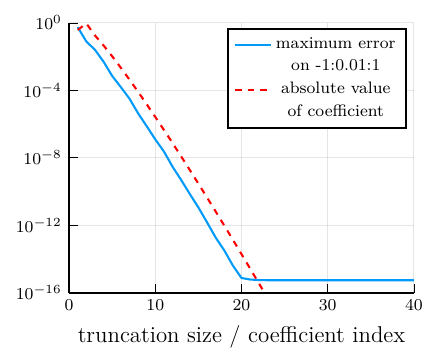}\label{FLO37}}
    \caption{Numerical results for \eqref{eq:exFIE2} with $\alpha=\beta=b=-1/2$ and $p = 3$.  
    }\label{FLO40}
\end{figure}

\end{example}

\begin{example}\label{exa3}
We consider the problem
\begin{equation}
u(x)+\lambda^2\mathcal{I}^{1/2}u(x)=1, \qquad x \in [-1,1], \label{eq:exFIE3}
\end{equation}
with the exact solution $u(x)=E_{1/2,1}(-\lambda^2\sqrt{1+x})$ (see \cref{thm6}), which can also be expressed as~\cite{haubold2011mittag}
\begin{equation}
u=  \exp(\lambda^4(1+x))\mathrm{erfc}(\lambda^2\sqrt{1+x}), \label{eq:uex3ex}
\end{equation}
where $\mathrm{erfc}$ denotes the complementary error function.
We shall compare the effect of $\lambda$ on the rate of convergence of the sum space method of~\cite{hale2018fast} and the JFP method. We let the constant $\lambda$ grow quadratically in \eqref{eq:exFIE3} since this is also the case in the time-fractional fractional heat/wave equation that we shall consider in \cref{ex:fheat}.

In the sum space method for \eqref{eq:exFIE3}, the solution is expanded as a sum of Legendre and weighted second-kind Chebyshev polynomials, or, expressed in terms of Jacobi polynomials, 
\begin{equation}
u(x) = \sum_{n=0}^{\infty} a_n P_n^{(0,0)}(x) + \sqrt{1+x}\sum_{n=0}^{\infty}b_nP_n^{(1/2,1/2)}(x) = \mathbf{S}(x) \bm{c}, \label{eq:changesum}
\end{equation}
where $\mathbf{S}$ is the `interleaved' sum space basis, i.e., 
\begin{equation}
\mathbf{S}(x) =\begingroup 
\setlength\arraycolsep{4pt} \left(\begin{array}{c | c | c | c | c}
P_0^{(0,0)}(x) & \sqrt{1+x}P_0^{(1/2,1/2)}(x) & P_1^{(0,0)}(x) & \sqrt{1+x}P_1^{(1/2,1/2)}(x) & \cdots
\end{array}    \right)\endgroup, \label{eq:Sbasis}
\end{equation}
and $\bm{c} = \left(
a_0 , b_0 , a_1 , b_1 , \cdots\right)^\top$ are the interleaved coefficients\footnote{We assume the coefficients $\bm{c}$ decay sufficiently fast so that the implicit change in the order of summation in \eqref{eq:changesum} is valid. }.   As shown in~\cite{hale2018fast}, the half-order integral operator $\mathcal{I}^{1/2}$ maps $\mathbf{S}$ to itself via a tridiagonal matrix, say $\widehat{I}_{1/2}$. Therefore, in the sum space basis, \eqref{eq:exFIE3} becomes 
$\mathbf{S}\left( \mathit{1} + \lambda^2 \widehat{I}_{1/2} \right)\bm{c} = \mathbf{S}\bm{e}_0$,
and the coefficients $\bm{c}$ can be obtained by solving the tridiagonal system 
\begin{equation}
\left( \mathit{1} + \lambda^2 \widehat{I}_{1/2} \right)\bm{c} = \bm{e}_0. \label{eq:sstridiag}
\end{equation}
 By contrast, with the JFP method, the system \eqref{eq:MLeqsyst} (with $\lambda \mapsto \lambda^2$)  has bandwidths $(k_*, N)$, where $N$ is the size of the truncated system and $k_* = \mu p = p/2$. 
 
\cref{FLO18} gives comparisons of the numerical results for solving \eqref{eq:exFIE3} with the sum space method and the JFP method (with $\alpha = \beta = 0$, $b = 0$ and $p = 2$ in \eqref{eq:MLeqsyst}).  

\begin{figure}[!htp]
    \centering
    \subfloat[Sum space method: double precision]{\includegraphics[width=0.45\textwidth]{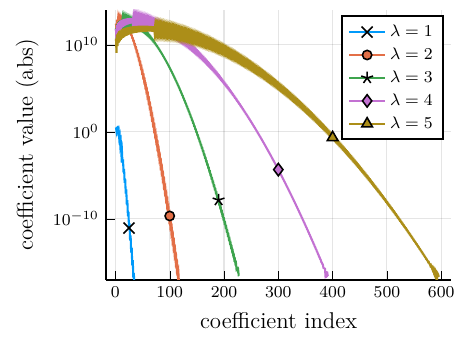}\label{FLO14}}
    \quad
    \subfloat[Sum space method: high precision]{\includegraphics[width=0.45\textwidth]{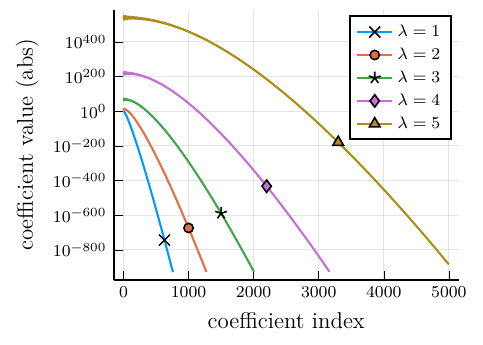}\label{FLO15}}\\
    \subfloat[JFP method]{\includegraphics[width=0.45\textwidth]{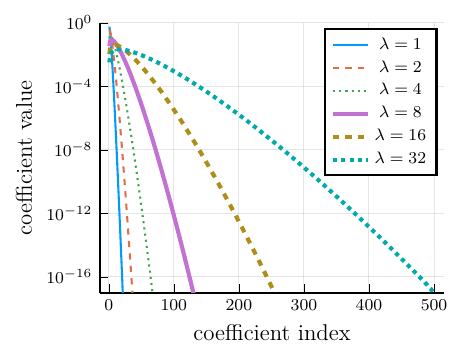}\label{FLO23}}
    \quad
    \subfloat[JFP method: maximum errors on the grid $-1$:0.01:1]{\includegraphics[width=0.45\textwidth]{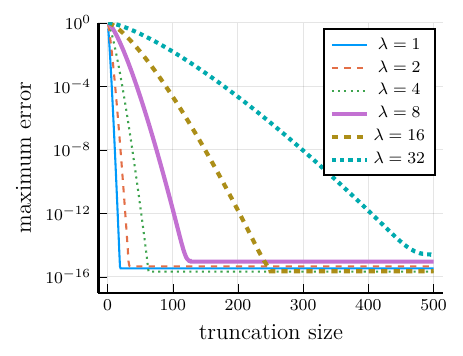}\label{FLO10}}\\
    \subfloat[Sum space method]{\includegraphics[width=0.45\textwidth]{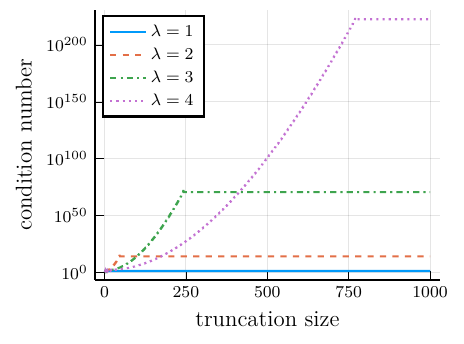}\label{FLO60}}
	\quad
	\subfloat[JFP method]{\includegraphics[width=0.45\textwidth]{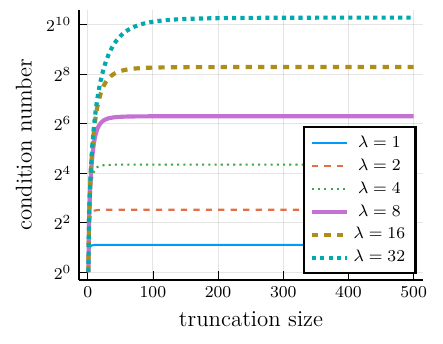}\label{FLO57}}
    \caption{
    	The top row shows  sum space coefficients obtained by solving \eqref{eq:sstridiag} with, \protect\subref{FLO14}: double precision  and \protect\subref{FLO15}: 3072-bit precision (see \cref{flo36}). 
    	The largest coefficient in \protect\subref{FLO15} for each $\lambda$ is approximately $\exp(2\lambda^4)$, which, for $\lambda = 2, \ldots, 5$, is on the order of $10^{14}$, $10^{70}$, $10^{222}$ and $10^{543}$, respectively. The JFP coefficients in \protect\subref{FLO23} are all bounded below $1$, even for much larger values of $\lambda$. \protect\subref{FLO10} shows the error of the JFP solution to \eqref{eq:exFIE3}, evaluated in double precision. For the JFP method, the truncation size (or number of coefficients) required to achieve $10^{-14}$ accuracy grows linearly in $\lambda$, whereas for the sum space method, the truncation size required for a fixed accuracy grows as $\mathcal{O}(\lambda^4)$.  For the sum space method \protect\subref{FLO60}, the ($2$-norm) condition numbers of \eqref{eq:sstridiag} grow roughly as $\exp(\mathcal{O}(n^2))$, where $n$ is the truncation size, and plateaus at a maximum value of approximately $\exp(2\lambda^4)$. For the JFP method \protect\subref{FLO57}, the condition numbers of \eqref{eq:MLeqsyst} asymptote to a value of approximately $1.8483\lambda^2$.
    }\label{FLO18}
\end{figure}

When using double precision in the sum space method, \cref{FLO13} shows that for $\lambda = 2, 3$, numerical errors in the solution on the order of $10^{-2}$ (and larger) are visible. This is because some coefficients are as large as $10^{15}$ (see \cref{FLO14}), which can lead to large cancellation errors. Using high precision, see \cref{FLO15}, shows how astronomically large the sum space coefficients become as $\lambda$ increases, which requires higher precision to compute and larger truncation sizes to ensure the solution achieves a specified accuracy. More precisely,  solving an $N\times N$ truncated version of the tridiagonal system \eqref{eq:sstridiag} in $q$-bit precision has a complexity of $\mathcal{O}(Nq\log q\log\log q)$. Computing an accurate solution to \eqref{eq:exFIE3} with the sum space method requires $q$-bit precision, with $q= \mathcal{O}(\lambda^4)$, and the required truncation size $N$ also grows as $\mathcal{O}(\lambda^4)$, hence the overall complexity of the sum space method is $\mathcal{O}(\lambda^8\log\lambda\log\log\lambda)$.

\begin{figure}[!htp]
    \centering
    \subfloat[Sum space method: numerical solutions]{\includegraphics[width=0.45\textwidth]{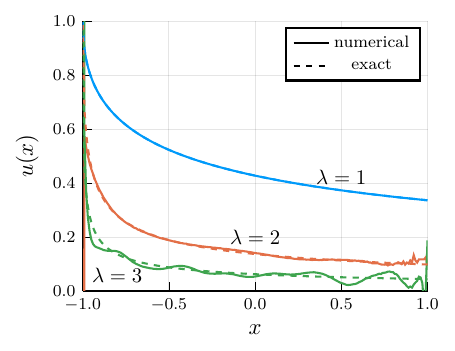}\label{FLO13}}
    \caption{The instability of sum space solutions of \eqref{eq:exFIE3}, computed in double precision, are clearly visible for $\lambda = 2, 3$.
   }\label{FLO16}
\end{figure}

In the JFP method for \eqref{eq:exFIE3},   $\mathcal{O}(N^2 q\log q \log\log q)$ operations are expended in computing the first $N$ columns of the matrix $I^{(0,0)}_{0,2,1/2}$ in $q$-bit precision using the pseudo-stabilized Algorithm 2 (see \cref{ps}). Solving the resulting $N \times N$ truncated version of the system \eqref{eq:MLeqsyst} for \eqref{eq:exFIE3} in $\widetilde{q}$-bit precision costs $\mathcal{O}(N^2 \widetilde{q}\log \widetilde{q} \log\log \widetilde{q})$ operations. As shown in \cref{ps}, $q = \mathcal{O}(N)$ and for \eqref{eq:exFIE3}, $N = \mathcal{O}(\lambda)$ and $\widetilde{q}$ is independent of $\lambda$ with\footnote{For \eqref{eq:exFIE3}, we solved the system \eqref{eq:MLeqsyst} in double precision, thus $\widetilde{q} = 53$.} $\widetilde{q}<q$,   hence the JFP method has an overall complexity of $\mathcal{O}(\lambda^3 \log \lambda \log\log \lambda)$, which is several orders lower than that of the sum space method ($\mathcal{O}(\lambda^8 \log \lambda \log\log \lambda)$ for \eqref{eq:exFIE3}). To reiterate, notwithstanding the fact that the JFP method requires the solution of a lower banded system (which, in addition requires high precision to compute), whereas the sum space method only requires the solution of an explicitly known tridiagonal system, the JFP method still has much lower complexity  than the sum space method for \eqref{eq:exFIE3}. 

The numerical results in \cref{FLO15} indicate that the magnitude of the largest sum space coefficient is similar to that of the largest power series coefficient in the normalized monomial basis, i.e., on the order of $\mathcal{O}\left(\exp(2\lambda^4)\right)$ (see \eqref{eq:maxmonoc} but with $\lambda$ replaced by $\lambda^2$ and $\mu = 1/2$). The normalized monomial and sum space coefficients of the solution to \eqref{eq:exFIE3} can be related using the following: we have that $\mathbf{P}^{(\alpha,\beta)} = \widetilde{\mathbf{M}}_0\widetilde{C}^{(\alpha,\beta)}$, where $\widetilde{C}^{(\alpha,\beta)}_{k,n} = 2^k C^{(\alpha,\beta)}_{k,n}$, see \labelcref{eq:normM0,eq6,eq:Ckn}.
Then since $u=E_{1/2,1}(-\lambda^2\sqrt{1+x})$, and using the sum space expansion \eqref{eq:changesum},
\begin{align*}
u(x) &= \sum_{n = 0}^{\infty} \frac{\left(-\sqrt{2}\,\lambda^2\right)^n}{\Gamma(1 + n/2)}\left(\frac{1+x}{2}\right)^{n/2} \\
&  = \sum_{n = 0}^{\infty} \frac{\left(2\lambda^4\right)^n}{\Gamma(1 + n)}\left(\frac{1+x}{2}\right)^{n} - \lambda^2 \sqrt{1+x} \sum_{n = 0}^{\infty} \frac{\left(2\lambda^4\right)^n}{\Gamma(3/2 + n)}\left(\frac{1+x}{2}\right)^{n} \\
& :=  \widetilde{\mathbf{M}}_0(x)\bm{c}_{\mathrm{e}} - \lambda^2 \sqrt{1+x} \,\widetilde{\mathbf{M}}_0(x)\bm{c}_{\mathrm{o}} \\
& = \mathbf{P}^{(0,0)}(x)\bm{a} + \sqrt{1+x}\,\mathbf{P}^{(1/2,1/2)}(x)\bm{b} \\
& = \widetilde{\mathbf{M}}_0(x)\widetilde{C}^{(0,0)}\bm{a} +  \sqrt{1+x}\,\widetilde{\mathbf{M}}_0(x)\widetilde{C}^{(1/2,1/2)}\bm{b}.
\end{align*}

Hence, $\bm{a} = \left(  \widetilde{C}^{(\alpha,\beta)}  \right)^{-1}\bm{c}_{\mathrm{e}}$ and $\bm{b} = -\lambda^2\left(  \widetilde{C}^{(\alpha,\beta)}  \right)^{-1} \bm{c}_{\mathrm{o}}$, where $\bm{a}$ and $\bm{b}$ are the sum space coefficients and the entries of $\bm{c}_{\mathrm{e}}$ and $\bm{c}_{\mathrm{o}}$ are $\bm{c}_{\mathrm{e}}(n) = \left(2\lambda^4\right)^n/\Gamma(1 + n)$ and $\bm{c}_{\mathrm{o}}(n) =   \left.\left(2\lambda^4\right)^n\middle/\Gamma(3/2 + n)\right.$. Since $\widetilde{\mathbf{M}}_0 = \mathbf{P}^{(\alpha,\beta)}\left(  \widetilde{C}^{(\alpha,\beta)}  \right)^{-1}$, the $n$-th column of $\left(  \widetilde{C}^{(\alpha,\beta)}  \right)^{-1}$ represents the Jacobi expansion coefficients of $((1+x)/2)^n$, which are bounded below one and strictly positive\footnote{Numerically, we found that $\|\left(  \widetilde{C}^{(\alpha,\beta)}_{0:n,0:n}  \right)^{-1}\|$ grows as $\mathcal{O}(\sqrt{n})$ for $\alpha = \beta = 0$ and as $\mathcal{O}(n^{1/4})$ for $\alpha = \beta = 1/2$.}. Since the entries of $\bm{c}_{\mathrm{e}}$ and $\bm{c}_{\mathrm{o}}$ are also strictly positive and their largest entries have a magnitude on the order of $\mathcal{O}\left(\exp(2\lambda^4)\right)$, we expect the largest sum space coefficients to have a similar magnitude, which is indeed the case in \cref{FLO15}.

The numerical results in \cref{FLO15} indicate that the truncation size required for the sum space method to achieve a prescribed accuracy grows as $\mathcal{O}(\lambda^4)$. By the same reasoning used above, we also expect this to be the case if the solution is expressed in the normalized monomial basis. Indeed, if we require $\vert \widetilde{\bm{u}}_{\nu,\mu}(n) \vert = \epsilon$ for large $\lambda$, with $\widetilde{\bm{u}}_{\nu,\mu}(n)$ defined in \eqref{eq:mutilde}, then using Stirling's approximation in the limit $n \to \infty$, we deduce that $n = \mathcal{O}(\lambda^{1/\mu})$, $\lambda \to \infty$. Replacing again $\lambda$ with $\lambda^2$ and setting $\mu = 1/2$, the truncation size $n$ grows as $\mathcal{O}(\lambda^4)$, as it does for the sum space method in \cref{FLO15}.   

\begin{remark}
We have considered the problem \eqref{eq:exFIE3} with order $1/2$, however similar conclusions hold if the problem has order $0 < \mu < 1$: the maximum sum space coefficient grows as $\mathcal{O}\left(\exp(2\lambda^{2/\mu})\right)$, and the required precision and truncation size scale as $\mathcal{O}(\lambda^{2/\mu})$, which leads to a complexity of $\mathcal{O}(\lambda^{4/\mu}\log\lambda\log\log\lambda)$. 
\end{remark}

For the JFP method, the magnitude of the coefficients in \cref{FLO23} and the $\mathcal{O}(\lambda)$ growth of the truncation size in \cref{FLO23,FLO10} can be derived from the bound \eqref{eq:chebcoeffsb} on the coefficients\footnote{The bound \eqref{eq:chebcoeffsb} is applicable to Chebyshev polynomials, however the JFP results in \cref{FLO18} were computed with Legendre polynomials (because the condition $\beta - b \in \mathbb{N}_0$ arising from Algorithm~2 rules out the use of Chebyshev polynomials; with Algorithm~1, however, Chebyshev polynomials are permissible). We use the bound \eqref{eq:chebcoeffsb} since it is simpler than the bounds on Legendre coefficients given in~\cite{zhao2013sharp}. Furthermore, the first $n$ Chebyshev and Legendre coefficients of the solution to \eqref{eq:exFIE3} can be related as follows, $\bm{u}_{0,0,0,2}(1\!\!\,:\,\!\!n) = R\bm{u}_{-1/2,-1/2,0,2}(1\!\!\,:\,\!\!n)$, where $R$ is an upper triangular matrix with~\cite{gautschi1972condition}  $\| R \| = \mathcal{O}(\sqrt{n})$. Hence, the Legendre coefficients are at worst larger than the Chebyshev coefficients by a factor of $\sqrt{n}$, which is inconsequential since the Chebyshev coefficients decay super-exponentially. }. Using \eqref{eq:uex3ex} and the relation $1+x = (1+y)^2/2$ from \eqref{eq31} with $p = 2$, it follows that $u$ attains its maximum modulus on the Bernstein ellipse $E_{\rho}$ in the complex $y$-plane at $y = -(\rho + \rho^{-1})/2$ (i.e., where $E_{\rho}$ intersects the negative real axis), hence in the bound \eqref{eq:chebcoeffsb},
\begin{equation}
M = \exp\left(\frac{\lambda^4}{2}\left(1   -(\rho + \rho^{-1})/2\right)^2  \right)\mathrm{erfc}\left( \frac{\lambda^2}{\sqrt{2}}\left(1 -(\rho + \rho^{-1})/2  \right)  \right).  \label{eq:Mexp}
\end{equation}
Setting $\rho = 1 + \epsilon$, $M$ becomes
\begin{equation}
M = \exp\left( \frac{\lambda^4 \epsilon^4}{8(1+\epsilon)^2}    \right)\left(1 + \mathrm{erf}\left( \frac{2\lambda^2\epsilon^2}{4\sqrt{2}(1 + \epsilon)}  \right)\right) \leq 2 \exp\left( \frac{\lambda^4 \epsilon^4}{8(1+\epsilon)^2}    \right).  \label{eq:Mbound}
\end{equation}
Hence $M \to 1$ as $\epsilon \to 0$ (or $\rho \to 1$) and we conclude from \eqref{eq:chebcoeffsb} that the JFP coefficients are bounded above by $2$ for any $\lambda$, which agrees with the numerical results in \cref{FLO23}. 

Next we estimate the growth of the truncation size as a function of $\lambda$ for the JFP method using the bounds \eqref{eq:chebcoeffsb} and \eqref{eq:Mbound}. First, we set $2M\rho^{-n} = c$, with $c \ll 1$ and for fixed $\lambda, n \gg 1$, estimate the value of $\rho$ that minimizes $2M\rho^{-n}$. 
We find that the requirement $\frac{\partial}{\partial\rho} M \rho^{-n} = 0$ leads to the following leading order estimate for the optimal value of $\rho$, 
\begin{equation*}
\rho = 1 + \epsilon \sim 1 + \left(\frac{2n}{\lambda^4}\right)^{1/3}, \qquad n, \lambda \gg 1.
\end{equation*}
Substituting this $\rho$ into $2M\rho^{-n} = c$, we find that
\begin{equation*}
n \sim \frac{\left(-\log(c/2)\right)^{3/4}}{2^{1/4}}\lambda, \qquad \lambda \to \infty,
\end{equation*}
hence the truncation size grows linearly in $\lambda$.

To make another comparison between the JFP and sum space bases, we take the perspective of orthogonal polynomials on algebraic curves~\cite{olver2021orthogonal,fasondini2021orthogonal,MF22}. In view of \eqref{eq:exFIE3}, in which $\mu = 1/2$ and with $p = 2$ in \eqref{eq31},  we consider orthogonal polynomials on the algebraic curve 
\begin{equation}
\gamma := \left\{ (x,y) \in \mathbb{R}^2 \: : \: \frac{1+x}{2} = \left(\frac{1+y}{2}\right)^2, \: x \in [-1, 1]\right\},
\end{equation}
which is shown in \cref{flo20}. We define the following inner product on $\gamma$,
\begin{equation*}
\langle f, g \rangle_{\gamma,w} := \int_{\gamma} f(x,y)g(x,y) w(x,y) d\sigma,
\end{equation*} 
where $d\sigma$ defines the arc length measure on the curve $\gamma$, hence we can set $d\sigma = \sqrt{1 + (x')^2}dy$, where $x' = \frac{dx}{dy} = y+1$. Setting $x = (1+y)^2/2 - 1$, defining $\widetilde{f}(y) = f(x,y)$, $\widetilde{g}(y) = g(x,y)$ and choosing $w(x,y)$ such that $w(x,y)d\sigma = w_{\alpha,\beta}(y)dy$, where the Jacobi weight is defined in \eqref{eq:jacw}, the inner product on the curve $\gamma$ becomes the standard Jacobi inner product in $y$, i.e., 
\begin{equation}
\langle f, g \rangle_{\gamma,w} = \int_{-1}^{1} \widetilde{f}(y)\widetilde{g}(y) w_{\alpha,\beta}(y) dy. \label{eq:jacipy}
\end{equation}
Hence, the JFP basis is a weighted orthogonal polynomial basis on the curve $\gamma$ (for \eqref{eq:exFIE3}, the JFP basis is not weighted because we set $b = 0$). The sum space basis \eqref{eq:Sbasis}, however, is not orthogonal with respect to  the inner product \eqref{eq:jacipy} for any choice of $\alpha, \beta > -1$. For example, setting $\alpha = \beta = 0$  and $1+x = (1+y)^2/2$ in \eqref{eq:jacipy}, we find that
\begin{equation*}
\langle P^{(0,0)}_n, P^{(0,0)}_m \rangle_{\gamma,w} = \int_{-1}^{1} P^{(0,0)}_n(x) P^{(0,0)}_m(x) \frac{dx}{\sqrt{2(1+x)}} \neq 0, \qquad n, m \geq 0,
\end{equation*}
similarly, 
\begin{equation*}
\langle P^{(0,0)}_n, \sqrt{1+\diamond}P^{(1/2,1/2)}_m \rangle_{\gamma,w} \neq 0,  \langle \sqrt{1+\diamond}P^{(1/2,1/2)}_n , \sqrt{1+\diamond}P^{(1/2,1/2)}_m \rangle_{\gamma,w} \neq 0.
\end{equation*}
The (non-orthogonal) sum space basis is an example of a frame~\cite[Example 3]{adcock2019frames}, which is a generalization of the notion of a basis.
 Finally, we note that the JFP and sum space bases for the problem \eqref{eq:exFIE3} are related as follows $\mathbf{S}(x) = \mathbf{P}^{(\alpha,\beta)}(y)R$, where $1+x = (1+y)^2/2$, hence 
 $R = \left(\| \mathbf{P}^{(\alpha,\beta)} \|^2\right)^{-1} \langle \mathbf{P}^{(\alpha,\beta)}, \mathbf{S}\rangle_{w_{\alpha,\beta}}$,
  which is an upper triangular matrix.

\begin{remark}\label{rem:pval}
    It is worth noting that if $k_*$, and thus $p$, are large, e.g., if $\mu = 1/2$ and $k_* = 5$ (and thus $p = 10$), then compared to choosing $k_*=1$, the rate of convergence is slower for small $\lambda$ but significantly faster for large $\lambda$. The effect of $k_*$ and $\lambda$ on the rate of convergence of the JFP method can likely be understood by using \eqref{eq:exmappedsol} and analyzing the bound  \eqref{eq:chebcoeffsb} as a function of $k_*$ and $\lambda$.
\end{remark}

\end{example}

\begin{example}\label{ex:fheat}
We consider a periodic one-dimensional time-fractional heat/wave equation of Caputo type,
\begin{equation}\label{eq32}
\left\{\begin{array}{l}
    \mathcal{D}^{\mu}_{\mathrm{C},t}u(x,t)=\mathcal{D}^2_xu(x,t),\quad -\infty<x<\infty, \quad 0<t<T, \quad 0 < \mu \leq 2,\\
    u(x,0)=f(x)=f(x+2\pi),
\end{array}\right.
\end{equation}
where $f(x)$ has a Fourier series $f(x)=\sum_{n=-\infty}^{+\infty}f_ne^{inx}$. For $0<\mu<1$, the equation is known as a fractional heat (or diffusion) equation. For $1<\mu<2$, \eqref{eq32} is referred to as a fractional wave equation and, in addition, $\partial_t u$  needs to be specified at $t = 0$. For simplicity, we shall set $\partial_t u(x,0) = 0$.

Using separation of variables, the solution has the form $u(x,t)=\sum_{n=-\infty}^{+\infty}f_n u_n(t)e^{inx}$, where the $u_n$ satisfy the FDEs
\begin{equation}
    \mathcal{D}_{\mathrm{C}}^{\mu}u_n(t)+n^2 u_n(t)=0, \qquad 0<t<T, \qquad u_n(0)=1, \qquad n \in \mathbb{Z}. \label{eq:unfdes}
\end{equation}
Recalling from \eqref{eq:RLCdefs} that $\mathcal{D}_{\mathrm{C}}^{\mu} = \mathcal{I}^{1-\mu}\mathcal{D}$ for $0 < \mu \leq 1$ and $\mathcal{D}_{\mathrm{C}}^{\mu} = \mathcal{I}^{2-\mu}\mathcal{D}^2$ for $1 < \mu \leq 2$, we apply $\mathcal{I}^{\mu}$ to the FDEs \eqref{eq:unfdes} to convert them to FIEs:
\begin{equation}
    u_n(t)+n^2 \mathcal{I}^{\mu} u_n(t)=1, \qquad 0<t<T, \qquad n \in \mathbb{Z}, \label{eq:unfies}
\end{equation}
where we have used the conditions $u_n(0) = 1$, $u_n'(0) = 0$ (for $1 < \mu \leq 2$), the latter of which originates from setting $\partial_t u(x,0) = 0$. Adapting \cref{thm6} to the interval $[0, T]$, it follows that $u_n(t) = E_{\mu,1}(-n^2t^{\mu})$. In practice, we consider \eqref{eq:unfies} for $0  \leq n \leq N$, where $N$ is the  truncation size of the Fourier series that is required to approximate $f$ to a given accuracy. Using a change of variables in \eqref{eq:unfies}, or using the exact solution, it can be verified that $u_n(t) = u_N\left( \left(\frac{n}{N}\right)^{2/\mu} t \right)$, hence we only need to solve \eqref{eq:unfies} for $n = N$, which reduces the computational cost by a factor of $N$. Next we map the FIE \eqref{eq:unfies} with $n = N$ to the unit interval $[-1, 1]$ by setting  $t = \frac{T}{2}(1+s)$ and $\widetilde{u}_N(s) = u_N\left(\frac{T}{2}(1+s)\right)$, then 
\begin{equation*}
    \widetilde{u}_N(s)+N^2\left(\frac{T}{2}\right)^\mu\mathcal{I}^{\mu} \widetilde{u}_N(s)=1, \qquad s \in [-1, 1].
\end{equation*}
This shows that in order to compute solutions to the fractional PDE \eqref{eq32}, we need to compute solutions to the FIE considered in the previous example, \eqref{eq:exFIE3},  in the large-$\lambda$ regime. 
 
As an example, in \cref{FLO43,FLO44} we let $u(x,0) = f(x) = e^{-\cos{2x}+1/2\sin x}-2\sin\sin x$, for which the truncation size of its Fourier series is $N = 29$ in double precision, and for all values of $\mu$, we let $p = 5$. 

\begin{figure}[!htp]
    \centering
    \subfloat[$\mu=2$]{\includegraphics[width=0.33\textwidth]{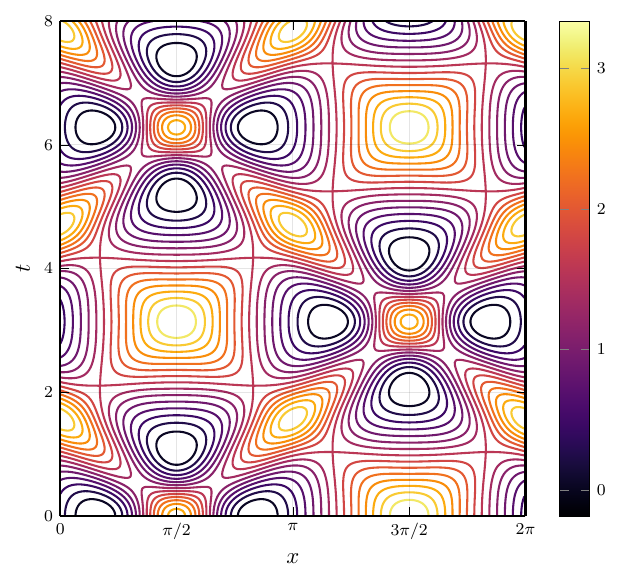}\label{FLO41}}
    \subfloat[$\mu=1.8$]{\includegraphics[width=0.33\textwidth]{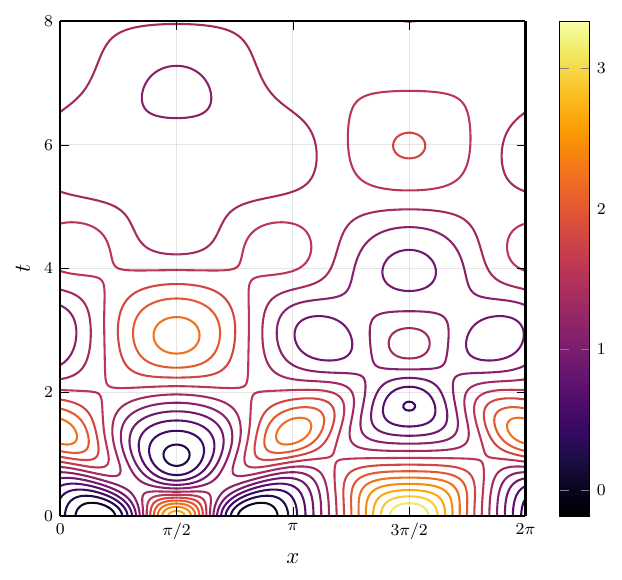}}
    \subfloat[$\mu=1.6$]{\includegraphics[width=0.33\textwidth]{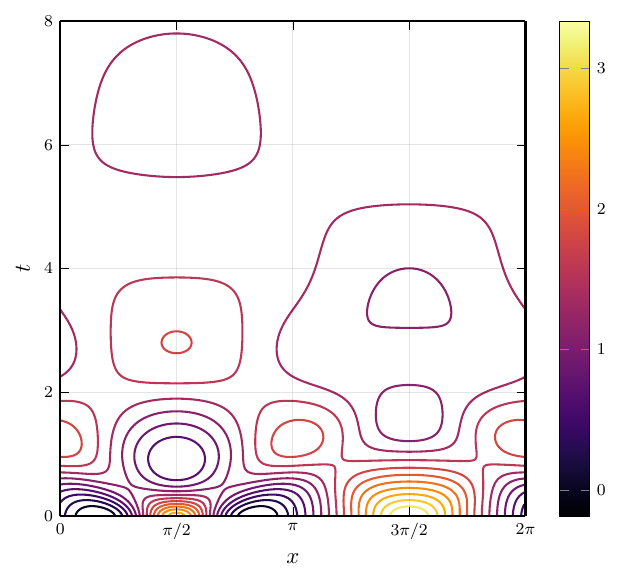}}
    
    \subfloat[$\mu=1.4$]{\includegraphics[width=0.33\textwidth]{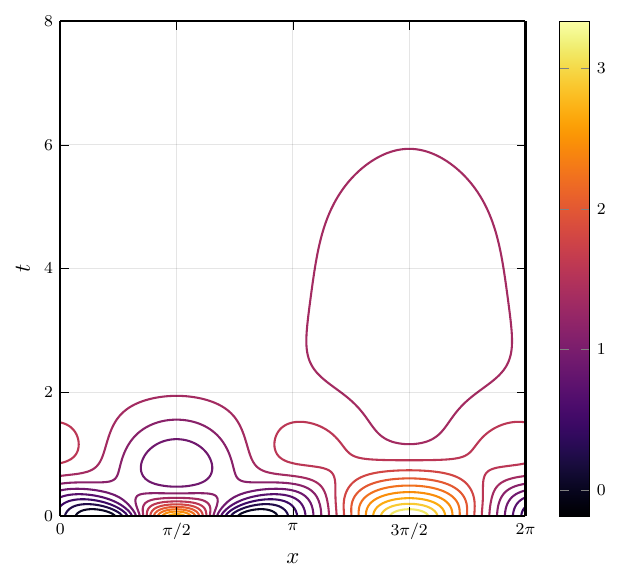}}
    \subfloat[$\mu=1.2$]{\includegraphics[width=0.33\textwidth]{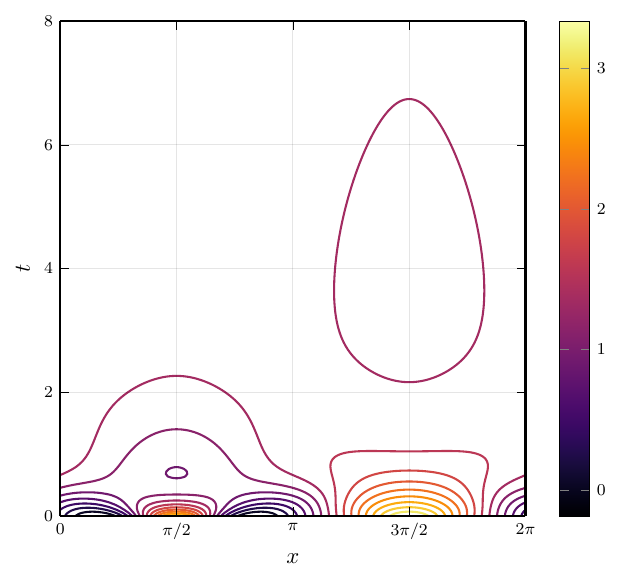}}
    \subfloat[$\mu=1$]{\includegraphics[width=0.33\textwidth]{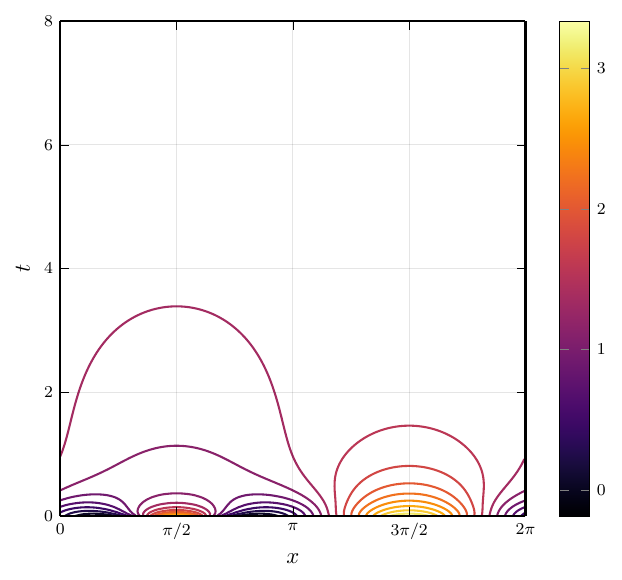}\label{FLO42}}
    \caption{Solutions to the time-fractional wave equation \eqref{eq32}, which show smooth transitions from the classical wave equation \protect\subref{FLO41} to the classical diffusion (or heat) equation \protect\subref{FLO42}.  
    }\label{FLO43}
\end{figure}

\begin{figure}[!htp]
    \centering
    \subfloat[$\mu=1$]{\includegraphics[width=0.33\textwidth]{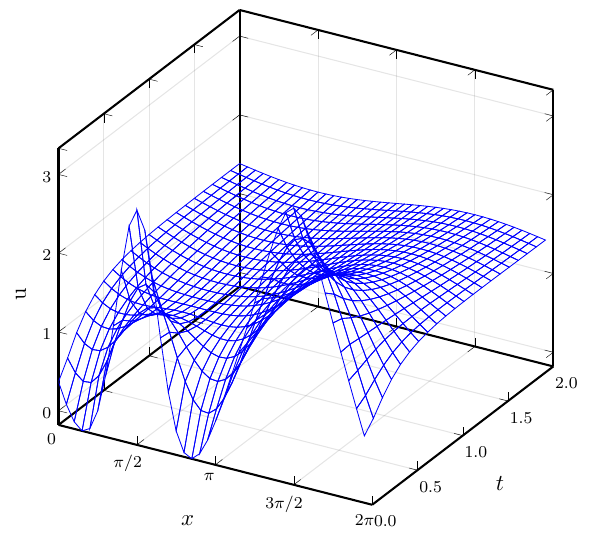}\label{classheat}}
    \subfloat[$\mu=0.8$]{\includegraphics[width=0.33\textwidth]{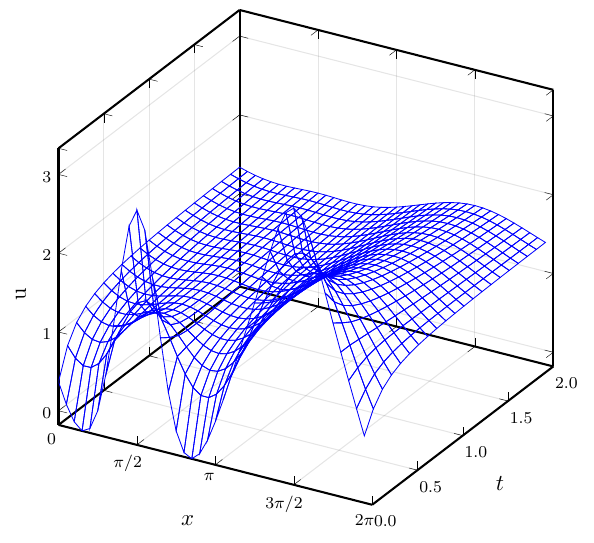}}
    \subfloat[$\mu=0.6$]{\includegraphics[width=0.33\textwidth]{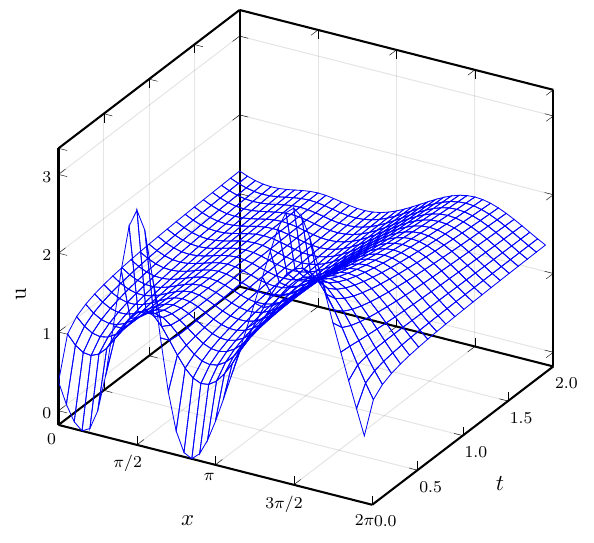}}
    
    \subfloat[$\mu=0.4$]{\includegraphics[width=0.33\textwidth]{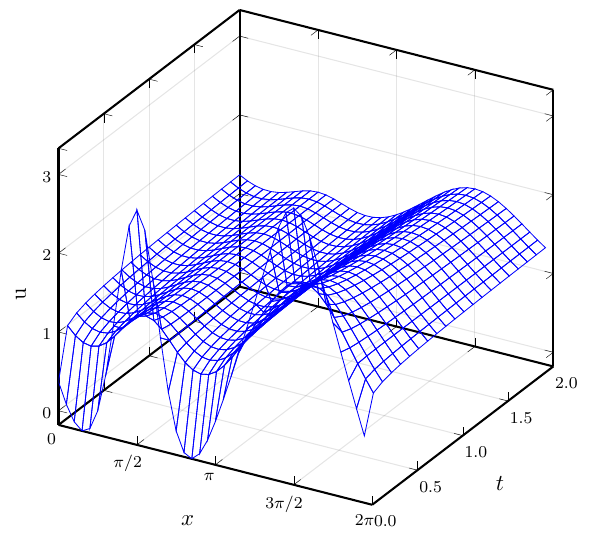}}
    \subfloat[$\mu=0.2$]{\includegraphics[width=0.33\textwidth]{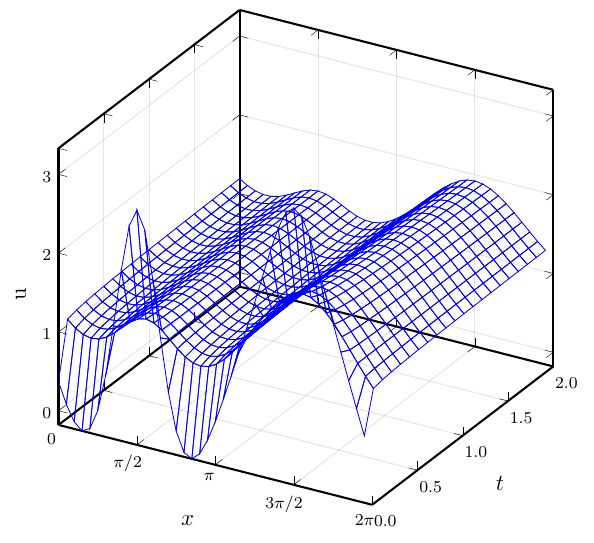}}
    \caption{Solutions to the time-fractional diffusion (or heat) equation \eqref{eq32}, including the classical diffusion equation \protect\subref{classheat}. As $\mu$ decreases from $1$, the  solutions decay faster initially but much more slowly for larger $t$. 
    }\label{FLO44}
\end{figure}
\end{example}

\subsection{More general FIEs and FDEs} \label{sect:probs}
In this subsection, we consider problems whose exact solutions are not known to us and include multiple orders (\cref{exa6}), variable coefficients (\cref{exa7}) and non-trivial boundary conditions (\cref{exa8}).

\begin{example}\label{exa6}
We consider a problem with multiple integer-order and fractional-order integral operators,
\begin{equation}
u(x)-\mathcal{I}^{1/2}u(x)+\mathcal{I}u(x)-\mathcal{I}^{3/2}u(x)+\mathcal{I}^{2}u(x)=1, \qquad x \in [-1,1], \label{eq:exmulto}
\end{equation}
which is also solved in \cite[Example 4]{hale2018fast}. From \eqref{eq4} it follows that the solution to this equation is an expansion in non-negative powers of $\sqrt{1 + x}$. Hence, \eqref{eq:urep} holds with $\mu = 1/2$, $\nu = 1$ and \eqref{eq:pbvals} gives the permissible choices of parameters. We let $p = 2$, $b = 0$ and $\alpha = \beta = 0$, i.e., we set $u = \mathbf{Q}^{(0,0,0,2)}\bm{u}_{0,0,0,2}$ in \eqref{eq:exmulto}, which leads to the lower banded system
\begin{equation}
\left(\mathit{1}-I^{(0,0)}_{0,2,1/2}+I^{(0,0)}_{0,2}-I^{(0,0)}_{0,2}I^{(0,0)}_{0,2,1/2}+\left(I^{(0,0)}_{0,2}\right)^2\right)\bm{u}_{0,0,0,2} = \bm{e}_0,  \label{eq:exmultosyst}
\end{equation} 
where we have made use of the semigroup property \eqref{eq20} to compute $\mathcal{I}^{3/2}$ and $\mathcal{I}^2$ in the JFP basis. Since $I^{(0,0)}_{0,2}$ has bandwidths $(2, 2)$, see \cref{thm4}, the matrix on the left-hand side of \eqref{eq:exmultosyst} has bandwidths $(4,\infty)$. 
The resulting solution (after truncating and solving \eqref{eq:exmultosyst})  and coefficients (indicating exponential convergence) are shown in \cref{FLO48}. 

\begin{figure}[!htp]
    \centering
    \subfloat[Solution]{\includegraphics[width=0.48\textwidth]{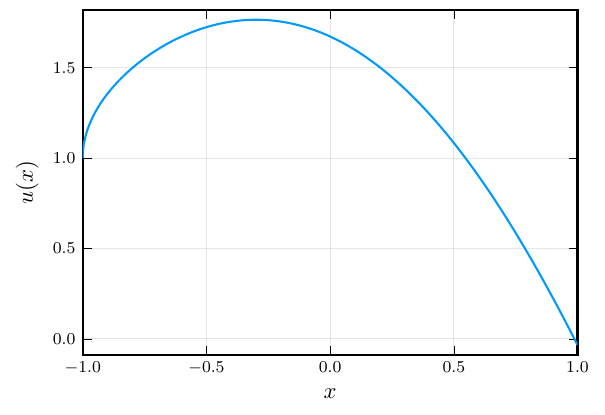}\label{FLO46}}
    \subfloat[Coefficients]{\includegraphics[width=0.48\textwidth]{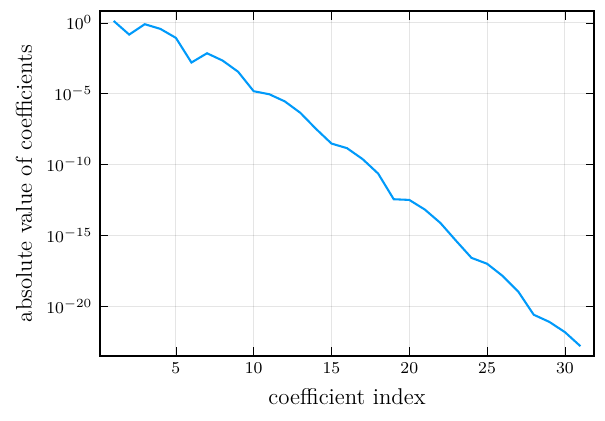}\label{FLO47}}
    \caption{ The solution to \eqref{eq:exmulto} in the JFP basis.
    }\label{FLO48}
\end{figure}

\end{example}

\begin{example}\label{exa7}
We consider an FIE with a variable coefficient,
\begin{equation} 
u(x)-\mathrm{erfc}\sqrt{1+x}\,\mathcal{I}^{1/2}u(x)=1, \qquad x \in [-1, 1],  \label{eq:exvcoeff}
\end{equation}
which is from \cite[Example 3]{hale2018fast}.  To express this FIE in the JFP basis, 
we require a matrix representing multiplication by a function $f(x)$, with $f(x) = \mathrm{erfc}\sqrt{1+x}$ in the case of \eqref{eq:exvcoeff}. Let  $\widetilde{f}(y)$ denote the function resulting from $f(x)$ after the change of variables defined in \eqref{eq31} and let $M_{\widetilde{f}}^{(\alpha,\beta)}$ be the matrix representing multiplication of the Jacobi  basis by $\widetilde{f}$ (hence, $\widetilde{f} \mathbf{P}^{(\alpha,\beta)} = \mathbf{P}^{(\alpha,\beta)}M_{\widetilde{f}}^{(\alpha,\beta)}$), then 
\begin{equation*}
\begin{split}
f(x)\mathbf{Q}^{(\alpha,\beta,b,p)}(x) &= (1+y)^b \widetilde{f}(y) \mathbf{P}^{(\alpha,\beta)}(y) = (1+y)^b  \mathbf{P}^{(\alpha,\beta)}(y)M_{\widetilde{f}}^{(\alpha,\beta)} \\
&= \mathbf{Q}^{(\alpha,\beta,b,p)}(x)M_{\widetilde{f}}^{(\alpha,\beta)}.
\end{split}
\end{equation*}
 The matrix $M_{\widetilde{f}}^{(\alpha,\beta)}$ can be constructed from a polynomial approximation to  $\widetilde{f}$ and Jacobi matrices in the manner described in \cite[section 2.3]{hale2018fast}. If $\widetilde{f}$ is approximated to a specified accuracy by a Jacobi polynomial expansion of degree $m$, then $M_{\widetilde{f}}^{(\alpha,\beta)}$ has bandwidths $(m, m)$. 

As in the previous example, the solution to \eqref{eq:exvcoeff} has an expansion in powers of $\sqrt{1 + x}$ and we again choose the parameters $b = \alpha = \beta = 0$, $p = 2$ for the JFP basis. Setting  $u = \mathbf{Q}^{(0,0,0,2)}\bm{u}_{0,0,0,2}$, \eqref{eq:exvcoeff} becomes
\begin{equation}
\left(  \mathit{1}- M_{\widetilde{f}}^{(0,0)}I^{(0,0)}_{0,2,1/2}   \right)\bm{u}_{0,0,0,2} = \bm{e}_0, \label{eq:exvcoeffsyst}
\end{equation}
where $\widetilde{f}(y) = \mathrm{erfc}((1+y)/\sqrt{2})$ is approximated on $[-1, 1]$ to machine epsilon in double precision with a degree $m = 21$ Legendre expansion. Hence, the matrix on the left-hand side of  \eqref{eq:exvcoeffsyst} has bandwidths $(m+1,\infty)$. The resulting exponentially convergent numerical solution and coefficients (after truncating and solving) are depicted in \cref{FLO51}. 

\begin{figure}[!htp]
    \centering
    \subfloat[Solution]{\includegraphics[width=0.48\textwidth]{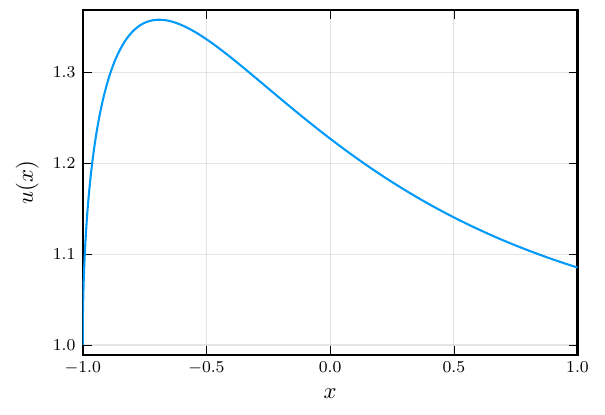}\label{FLO49}}
    \subfloat[Coefficients]{\includegraphics[width=0.48\textwidth]{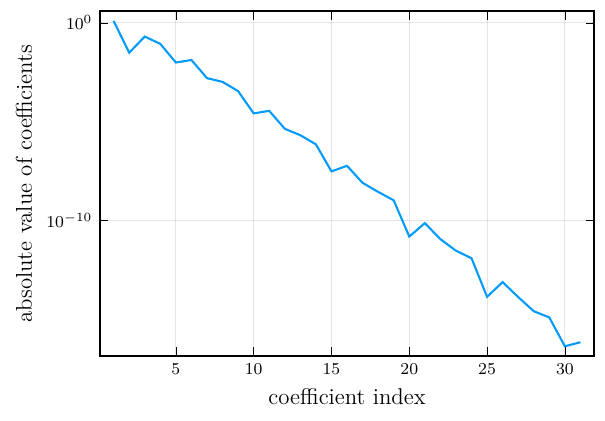}\label{FLO50}}
    \caption{The solution to \eqref{eq:exvcoeff} in the JFP basis.
    }\label{FLO51}
\end{figure}

\end{example}

\begin{example}\label{exa8}
The final problem we consider in this section is the classical Bagley--Torvik equation~\cite{bagley1983fractional,torvik1984appearance},
\begin{equation}
u''+\mathcal{D}^{1/2}u+u=0\quad\text{on}\quad(-1,1), \qquad u(-1)=1, u(1)=0,  \label{eq:exFDEBT}
\end{equation}
where $\mathcal{D}^{1/2}$ can be of either Riemann--Liouville (RL) or Caputo (C) type, defined in \eqref{eq:RLCdefs}. To reformulate \eqref{eq:exFDEBT} as an FIE, we let 
\begin{equation}
u(x)=\mathcal{I}^2 v(x)+a(1+x)+1,  \label{eq:uiiv}
\end{equation}
hence one boundary condition is satsfied, $u(-1) = 1$, and $u''=v$. Noting that $\mathcal{D}^{1/2}_{\mathrm{RL}}u(x)=\mathcal{D}^{1/2}_{\mathrm{C}}u(x)+\frac{1}{\Gamma(1/2)\sqrt{1+x}}$, where $\mathcal{D}^{1/2}_{\mathrm{C}}u(x)=\mathcal{I}^{3/2}v(x)+\frac{a}{\Gamma(3/2)}\sqrt{1+x}$, we let
$$\mathcal{D}^{1/2}u(x)=\mathcal{I}^{3/2}v(x)+\frac{a}{\Gamma(3/2)}\sqrt{1+x}\left(+\frac{1}{\Gamma(1/2)\sqrt{1+x}}\right)_{\mathrm{RL}},$$ then \eqref{eq:exFDEBT} becomes the FIE
\begin{equation}
    (1+\mathcal{I}^{3/2}+\mathcal{I}^2)v(x)+\frac{a}{\Gamma(3/2)}\sqrt{1+x}+1\left(+\frac{1}{\Gamma(1/2)\sqrt{1+x}}\right)_{\mathrm{RL}}=0,  \label{eq:exFIEBT}
\end{equation}
with $a$ to be determined and $v$ subject to the boundary condition
\begin{equation}
\mathcal{I}^2 v(1)+2a+1=0. \label{eq:exFIEBTbcs}
\end{equation}

The solution to \eqref{eq:exFIEBT} has an expansion in $\{(1+x)^{n/2}\}_{n=-1}^{\infty}$, hence \eqref{eq:urep} holds with $\mu = 1/2$ and $\nu = 1/2$ and the permissible parameter values for the JFP basis are $p = 2 k_*$, $k_* =\mathbb{N}_+$ and $b  = -k_* -n$, $n \in \mathbb{N}_0$, with $b > -2k_*$ see \eqref{eq:pbvals}. We set $k_* = 1$, in which case $p = 2$, $b = -1$ and we let $\alpha = \beta = 0$. For brevity we omit the parameter values of the JFP basis $\mathbf{Q}^{(0,0,-1,2)}$ and set $v = \mathbf{Q}\bm{v}$, $f = \mathbf{Q}\bm{f}$ and $g = \mathbf{Q}\bm{g}$, where $f(x)=1\left(+\frac{1}{\Gamma(1/2)\sqrt{1+x}}\right)_{\mathrm{RL}}$ and $g(x)=(1+x)+\frac{\sqrt{1+x}}{\Gamma(3/2)}$. Substituting this into \eqref{eq:exFIEBT} and \eqref{eq:exFIEBTbcs}, we obtain 
$$\begingroup 
\setlength\arraycolsep{8pt}\left(\begin{array}{c |c}
 \mathbf{Q}(1)\left(I_{-1,2}^{(0,0)}\right)^2 & 2 \\
 \hline
    \mathit{1}+I_{-1,2,3/2}^{(0,0)}+\left(I_{-1,2}^{(0,0)}\right)^2 & \bm{g}
   
\end{array}\right)\endgroup 
\left(\begin{array}{c}
    \bm{v}\\
    a
\end{array}\right)
=
\left(\begin{array}{c}
    -1\\
    -\bm{f}
\end{array}\right),$$
where $\mathbf{Q}(1)=\left(\begin{array}{c c  c } Q_0(1) & Q_1(1) & \cdots\end{array}\right) \in \mathbb{R}^{1 \times \mathbb{N}_0}$. We truncate and solve this system and obtain an approximation to the solution of the original FDE from the relationship \eqref{eq:uiiv}. The exponentially convergent numerical results are given in \cref{FLO56} and agree with those shown in \cite[Examples 7 and 9]{hale2018fast}.

\begin{figure}[!htp]
    \centering
    \subfloat[Solution]{\includegraphics[width=0.48\textwidth]{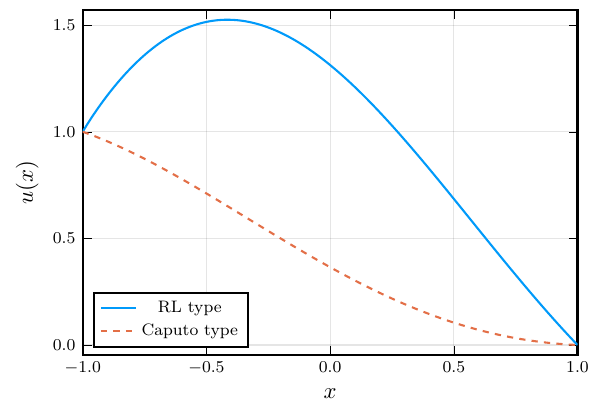}\label{FLO54}}
    \subfloat[Coefficients of $\mathcal{I}^2v$]{\includegraphics[width=0.48\textwidth]{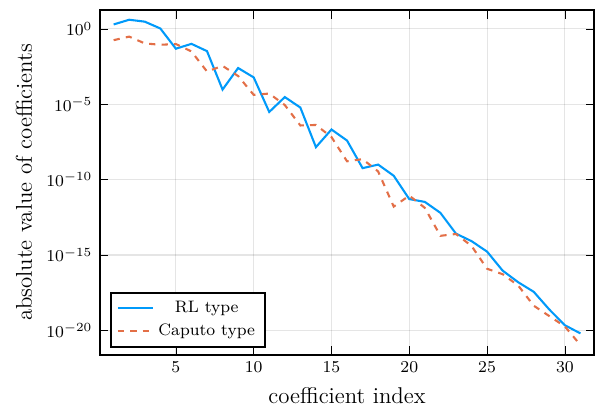}\label{FLO55}}
    \caption{The solution to the Bagley--Torvik FDE \eqref{eq:exFDEBT} by reformulation as an FIE \eqref{eq:exFIEBT} and solved via the  JFP method.
    }\label{FLO56}
\end{figure}

\end{example}

\begin{remark}
The general form of the linear FIEs that can be solved with the JFP method is given by
\begin{equation}
a_0(x)u(x) + a_1(x)\mathcal{I}^{\mu_1}\left[b_1u \right](x) + \cdots +  a_N(x)\mathcal{I}^{\mu_N}\left[b_Nu  \right](x) = f(x), \label{eq:genFIE}
\end{equation}
where it is assumed that $\mu_1 > 0$ and $\mu_k$, $2 \leq k \leq N$ are positive rational multiples of $\mu_1$, hence $\mu_k = \frac{m_k}{n_k}\mu_1$, where the $m_k$ and $n_k$ are positive (and relatively prime) integers. We can set $\mu_1 = k_*/p$, where $k_* \in \mathbb{N}_+$ and $p> 0$ and we further assume that there exists a positive integer $n$ such that the given functions $a_0(x)$, $f(x)$, $a_k(x)$, $b_k(x)$, $1 \leq k \leq N$ have convergent expansions in non-negative powers of $(1+x)^{1/(np)}$, where $n \geq \max\{n_2, \ldots, n_N\}$. We allow $f(x)$ to have an expansion in the basis\footnote{\eqref{eq:exFIE2} is an example of this case with $\mu_1 = 1/3$ (we can choose $k_* = 1$ and $p=3$, for example) and $\nu = 3/2$.}  $\mathbf{M}_{\nu-1}^{1/(np)}$, $\nu > 0$, in which case we divide by $(1+x)^{\nu - 1}$ in \eqref{eq:genFIE} and we relabel $(1+x)^{1-\nu}u(x)$ as $u(x)$. Hence, we assume without loss of generality that $f$ has an expansion in non-negative powers of $(1+x)^{1/(np)}$.   Then the solution to \eqref{eq:genFIE} has an expansion in the JFP basis $\mathbf{Q}^{(\alpha, \beta,0,np)}$.
\end{remark}

\begin{remark}
We have demonstrated that the sum space method for the FIE \eqref{eq:exFIE3} performs poorly compared to the JFP method as $\lambda$ increases. However, as illustrated in \cite{hale2018fast}, the sum space method converges exponentially fast in linear complexity (because it yields banded or almost-banded systems) and in double precision for a wide range of problems (examples of which include \eqref{eq:exmulto} and \eqref{eq:exFDEBT}). For these problems the sum space method is preferable to the JFP method because it converges at a similar rate but with lower complexity (because, by comparison, the JFP method always leads to a lower banded system for fractional order problems) and without the need to compute integration matrices in high precision. The sum space method has linear complexity if the variable coefficients are analytic functions of $x$ but for an FIE such as \eqref{eq:exvcoeff}, the  method also leads to a lower banded system. 

The sum space method is only applicable to equations of rational order $\mu = p/q$ and requires the direct sum of $2q$ different weighted orthogonal polynomial bases\footnote{$q$ bases for the domain of the FDE/FIE and $q$ bases for the range, however the case $q = 2$ is an exception for which only $2$ bases are required, see \eqref{eq:Sbasis}}. For example, for \eqref{eq:exFIE2}, $q = 6$ because the solution has an expansion in powers of $(1 +x)^{1/6}$, and thus $12$ weighted bases are required while the JFP method uses a single basis (but with different parameter values $p$ and $b$)  for all fractional order (including irrational order) problems.

We have not attempted to make a rigorous classification of the types of problems for which the JFP method is superior to the sum space method and vice versa. However, \cref{exa3} suggests that the sum space method performs poorly for problems in which the largest monomial coefficient of the solution becomes large (e.g., on the order of $10^{10}$). Otherwise, we conjecture, the sum space method performs well.

Large condition numbers of the systems that arise in the sum space method are not good indicators of the accuracy of the resulting solutions. For example, in \cite[Example 11]{hale2018fast}, the following fractional Airy equation is considered, 
\begin{equation}
\epsilon i^{3/2} \mathcal{D}_{\mathrm{RL}}^{3/2} u(x) - x u(x) = 0, \qquad x \in (-1, 1), \qquad u(-1)=0, u(1)=1,  \label{eq:fractairy}
\end{equation}
which is a singularly perturbed problem with an increasingly oscillatory solution as $\epsilon \to 0$. For this problem the sum space method achieves high accuracy (around $10^{-10}$  for $\epsilon = 10^{-4}$) despite large condition numbers (on the order of $10^{15}$ and larger). This is reminiscent of the high accuracy achieved by the ultraspherical spectral method in \cite{olver2013fast} for a singularly perturbed standard Airy equation (i.e.,  with $i^{3/2} \mathcal{D}_{\mathrm{RL}}^{3/2}$ in \eqref{eq:fractairy} replaced by $\mathcal{D}^2$) despite the large condition numbers that arise. By our conjecture, we expect the sum space method to perform well for the problem \eqref{eq:fractairy}, large condition numbers notwithstanding, because the monomial coefficients of the solution to \eqref{eq:fractairy} do not become large. 
\end{remark}

%% file: conclusion.tex
\section{Conclusion}

We have illustrated the application of the JFP method to a variety of FIEs (including FDEs and a fractional PDE reformulated as FIEs) in which  exponentially fast convergence to the solution is achieved. The JFP method converges much faster and with a lower overall complexity than the sparse sum space method in~\cite{hale2018fast} for solutions whose power series (or shifted monomial) coefficients become large. For such problems, a relatively large number of coefficients are needed to resolve the solution in which case  the use of high-precision arithmetic is essential to scale up the JFP method while retaining sufficient accuracy. Pseudo-stabilization of the unstable algorithms we have introduced for fractional integration matrices incorporates high-precision  automatically (see~\cref{ps}) and increases the complexity of the JFP method from $\mathcal{O}(N^2)$ (without high-precision, but unstable)  to  $\mathcal{O}(N^3\log N\log\log N)$.  

\subsection{Future work and open problems}

We shall construct differentiation matrices for the JFP basis, which will be banded and could be used to solve FDEs and fractional PDEs  without having to reformulate them as FIEs, as we have done  in this paper.  Adding Newton iteration in function space~\cite{birkisson2012automatic}, the JFP method would also be applicable to nonlinear FIEs and FDEs.

The superior performance of the JFP method in \cref{exa3} of \cref{sect:testprobs} suggests it could be an effective method for computing  Mittag--Leffler functions. Just as the ultraspherical spectral method was shown in~\cite{crespo2020multidomain} to be an effective method for the global computation of a special function (the Gauss hypergeometric  function), the JFP method could be adapted in a similar way to efficiently evaluate Mittag--Leffler functions on intervals and regions in the complex plane.

The  convergence and stability analysis of the JFP method for general linear FIEs \eqref{eq:genFIE} is a topic for future research. This analysis could prove the bounds on condition numbers that we found in \cref{FLO57} and reveal the dependence of the rate of convergence on the parameter $p$ in the JFP basis (see \cref{rem:pval}). 

Another topic for future research is stable algorithms for fractional integration matrices in the JFP basis with optimal (i.e., $\mathcal{O}(N^2)$) complexity. As mentioned in \cref{sect:alg2sect}, one possibility is to use asymptotic approximations to the entries of the fractional integration matrices to stabilize the Sylvester equations in Algorithm~2.

%% file: ps.tex
\section{Pseudo-stabilization}\label{ps}

\subsection{Introduction}

Pseudo-stabilization is a technique that uses high-precision computations to ensure that the output of an unstable algorithm has a specified accuracy. Pseudo-stabilization is called thus because it makes an unstable algorithm appear stable by outputting accurate results. This technique is algorithm-specific and here we illustrate it for Algorithm~2 (see \cref{sect:alg2sect}), which computes the entries of fractional integration matrices via the recurrence relation \eqref{eq19} (with $A = I^{(\alpha,\beta)}_{b,p,\mu}$) in which errors grow exponentially. Since high-precision computations are more expensive than standard double precision, our aim is to use as low a precision as possible while ensuring the integration matrices are accurate to the specified tolerance. To choose this optimal precision, we need to know how fast errors grow in \eqref{eq19}, examples of which are shown in \eqref{flo30}. Since the exact values of the entries of $ I^{(\alpha,\beta)}_{b,p,\mu}$ are not known, throughout this appendix, errors in a given precision are computed with reference to results computed in a much higher precision. 

\subsection{Growth of errors}

\begin{figure}
    \centering
    \subfloat[Errors in each column]{\includegraphics[width=0.45\textwidth]{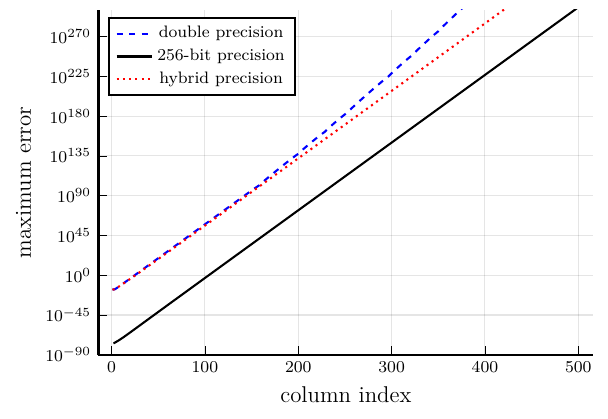}\label{flo28}}
    \quad
    \subfloat[Column-wise growth rate of errors]{\includegraphics[width=0.45\textwidth]{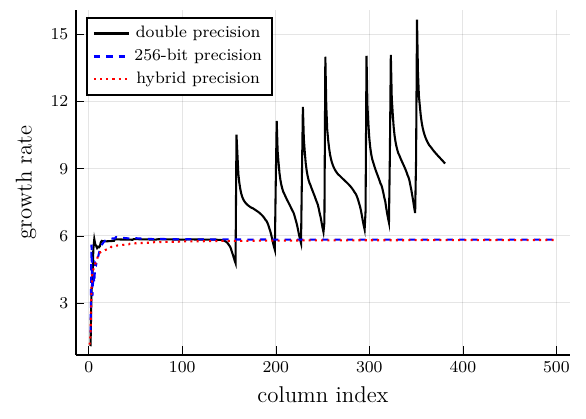}\label{flo29}}
    \caption{The errors from computing $I^{(0,0)}_{0,2,1/2}$   using Algorithm~2. For the hybrid precision results, we computed the first $p=2$ columns of $I^{(0,0)}_{0,2,1/2}$ in double precision (using Algorithm~1), while the recurrence coefficients (i.e., the entries of $B$ and $C$ in \eqref{eq19}) were computed with $256$-bit precision (see \cref{flo36}).  In \cref{flo29}, the growth rate is the  maximum error in a column divided by that in the previous column.  
 }\label{flo30}
\end{figure}

\begin{figure}[!htp]
    \centering
    \subfloat[$r(q,n)$ for $I_{0,3,1/3}^{(0,0)}$]{\includegraphics[width=0.33\textwidth]{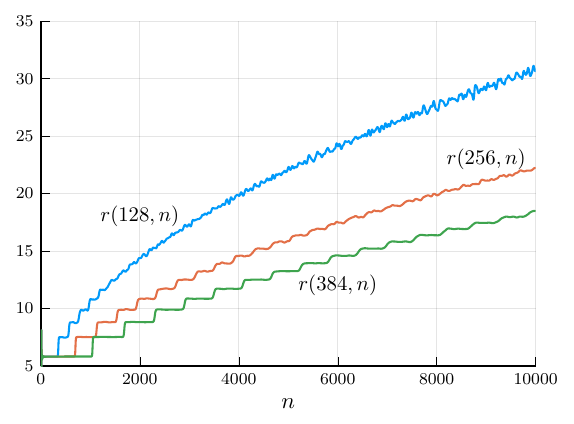}\label{flo32}}
    \subfloat[$r(q,n)$ for $I_{0,3,2/3}^{(0,0)}$]{\includegraphics[width=0.33\textwidth]{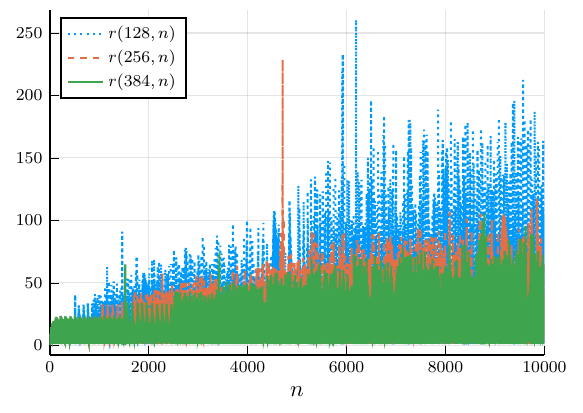}\label{flo33}}
    \subfloat[$\left(\frac{E(q,n)}{E(q,n-128)}\right)^{1/128}$ for $I_{0,3,2/3}^{(0,0)}$]{\includegraphics[width=0.33\textwidth]{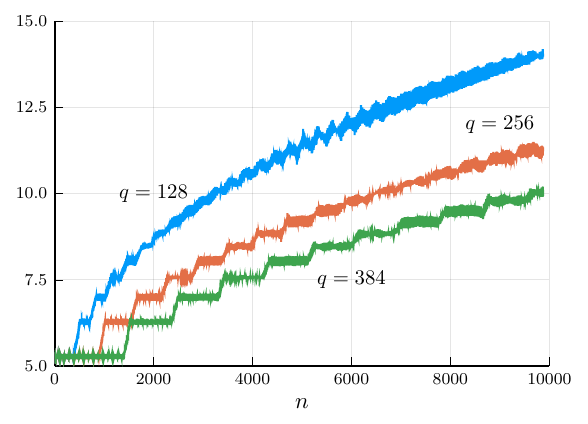}\label{flo35}}
    \caption{\protect\subref{flo32}: The estimated growth rate $r$ of the error (see \eqref{eq:rdef})  for the matrix $I_{0,3,1/3}^{(0,0)}$ resembles a superposition of step functions with the scaling property \eqref{eq26}.  \protect\subref{flo33}: $r$ is noisy for $I_{0,3,2/3}^{(0,0)}$ but the averaged version of $r$ in \protect\subref{flo35}, denoted by $r_0$, also satisfies \eqref{eq26}. 
    }\label{flo34}
\end{figure}

The error in the first $p$ columns of $I^{(\alpha,\beta)}_{b,p,\mu}$ is determined by the precision used in Algorithm~1, while the growth rate of the error in subsequent columns depends on the recurrence coefficients (i.e., the entries of $B$ and $C$ in \eqref{eq19}) and hence also on the precision with which they are computed, as confirmed by \cref{flo30}. For recurrence coefficients computed in double precision, \cref{flo29} shows that the errors grow roughly at a constant rate for the first circa $150$ columns, after which the growth rate starts shooting up. The errors in higher (256-bit) precision grow at roughly a constant rate in \cref{flo29}, however we shall find in \cref{flo34} that if a sufficiently large number of columns are computed, the growth rate of the error also increases (for any precision). Due to the complicated behavior of the growth rate of the error, we shall adopt an empirical approach to model and predict the growth of the error. 

\begin{figure}[!htp]
    \centering
    \includegraphics[width=0.6\textwidth]{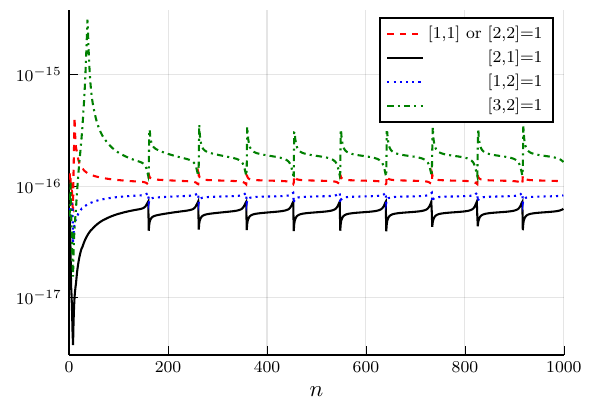}
    \caption{A plot of $E(q,n;e_0)/\widetilde{E}(q,n)$ in double precision (hence, $q = 53$ and $e_0 \approx 10^{-16}$) for the integration matrix $I_{0,2,1/2}^{(0,0)}$, which  illustrates the approximation \eqref{eq:simerr}. $\widetilde{E}(q,n)$ is not unique since there are multiple initial values which can be set to $0$ or $1$; the legend indicates the entry of $A$ in \eqref{eq19} that is set to $1$ to compute  $\widetilde{E}(q,n)$. 
  }
    \label{flo31}
\end{figure}

Let $q$ be the precision and $\epsilon$ the corresponding machine epsilon, then $\epsilon(q)=2^{1-q}$. We let $e_0$ and $E(q,n;e_0)$ denote, respectively, the maximum errors in the initial $p$ columns and in the $n$-th column of $I^{(\alpha,\beta)}_{b,p,\mu}$, with $n \geq 0$. \cref{flo28} shows examples of $E(q,n;e_0)$ for $q = 53$ and $q = 256$.  We let $E(q,n) = E(q,n;e_0)/e_0$, hence, by definition, $E(q,p-1) = 1$ and therefore we simulate the normalized maximum error $E(q,n)$ by $\widetilde{E}(q,n)$, which is defined as the maximum value in the $n$-th column of $A$ when one or more of the initial values in the recurrence \eqref{eq19} are set to $1$ and  the rest are set to zero.
We have found that 
\begin{equation}
\frac{E(q,n;e_0)}{e_0} = E(q,n) \approx   \widetilde{E}(q,n)   \label{eq:simerr}
\end{equation}  
holds, see \cref{flo31} for an example. This implies that we can model the growth of errors with $\widetilde{E}(q,n)$.  

Next we consider the estimated growth rate of the error,
\begin{equation}
\frac{E(q,n;e_0)}{E(q,n-1;e_0)} \approx \frac{\widetilde{E}(q,n)}{\widetilde{E}(q,n-1)} := r(q,n). \label{eq:rdef}
\end{equation}
\cref{flo29} shows $r(q,n)$ for $q = 53, 256$ and \cref{flo34} shows $r(q,n)$ for much larger $n$ (up to $10^4$). In \cref{flo29}, $r(256,n)$ is approximately constant, however in \cref{flo34} it is seen that $r(256,n)$ also increases in a step-like manner for large enough $n$. \cref{flo33} shows that $r$ can be very noisy, however the denoised or averaged version of $r$, which we denote by $r_0(q,n)$, takes the form of a superposition of step functions satisfying the scaling property (see \cref{flo32,flo35})
\begin{equation}\label{eq26}
r_0(q,n)\approx r_0(\lambda q,\lambda n).
\end{equation}
Taking $\lambda = 2$ as an example, we have 
\begin{align}
\widetilde{E}(2 q,2 n) &\approx r_0(2q,2n)\widetilde{E}(2 q,2 n-1) \approx r_0(2q,2n)r_0(2q,2n-1)\cdots r_0(2q,1) \label{eq:rtoE} \\
& \approx r_0^2(q,n)\cdots r_0^2(q,1) = \widetilde{E}^2(q,n),\notag
\end{align}
where we used the fact that $\widetilde{E}( q,0) = 1$. By \eqref{eq:simerr}, this implies that  $E(2 q,2 n) \approx E^2(q,n)$. In general, we use the  approximation
\begin{equation}\label{eq27}
E(\lambda q,\lambda n) \approx E^\lambda(q,n),
\end{equation}
which means that we can use a simulation on lower precision and fewer iterations to predict error propagation on higher precisions through a larger number of iterations.  The approximation \eqref{eq27} is exact for the simplest case in which $r(q,n)$ and $r(\lambda q,\lambda n)$ are constant and equal (e.g., for double precision with $q = 53$, $n \leq 150$, see  \cref{flo29}, and $\lambda > 1$). If we let $e_0(q) = 2^{-q}$, then for constant $r$, $E(n,q;e_0) = r^n2^{-q}$ and $E(\lambda n, \lambda q;e_0) = r^{\lambda n}2^{-\lambda q} = E^\lambda (n,q;e_0)$.

\begin{remark}
 In \cref{flo34}, the height of the steps of $r_0$ for different $q$ are the same but the lengths are proportional to $q$ (hence the property \eqref{eq26}). Except for the first step, which is the longest for every $q$, every step has the same length.
We are surprised by the step-like behavior of the growth rate.   It seems to be caused by the errors in computing the integration matrix $I_{b,p}^{(\alpha,\beta)}$. We conjecture that if we could compute the integration matrix exactly, then the first step of $r(q,n)$ would be infinitely long, i.e., the growth rate would be constant, which is what we would expect for a linear recurrence problem.
\end{remark}

\begin{remark}
We have found that even though the approximations \eqref{eq26} and \eqref{eq27} are rather crude, especially when $r$ is noisy, there is enough `noise cancellation' in the products of growth rates in \eqref{eq:rtoE}  (with $r_0$ replaced by the actual growth rates $r$ defined in \eqref{eq:rdef}) for the approximation \eqref{eq27} to yield useful estimates. 
\end{remark}

\subsection{Implementation}

We now use the empirical estimate \eqref{eq27} for the growth of errors to estimate the lowest precision required to meet a given tolerance. These estimates were used to compute the fractional integration matrices that appear in the problems in \cref{spec}.

Suppose we aim to achieve a maximum error of at most $\delta$ in the computation of the first $N+1$ columns of $A = I^{(\alpha,\beta)}_{b,p,\mu}$ using the recurrence \eqref{eq19}. Setting $e_0 = 2^{-q}$ and using \eqref{eq:rdef}, we seek a $q$ such that
\begin{equation}
E(q,N;e_0) = 2^{-q}E(q,N) < \delta.  \label{eq:qrel}
\end{equation}
Now the idea is to simulate $E$ with $\widetilde{E}$ on a chosen lower precision $q_0$ and a smaller number of iterations $m < N$ (which is to be determined) and use  \eqref{eq27} to estimate the smallest $q$ satisfying \eqref{eq:qrel}.

Using \eqref{eq27} with $\lambda = q/q_0$, we have $E(q,N) \approx E(q_0,Nq_0/q)^{q/q_0}$, hence we require 
\begin{equation*}
E(q_0,Nq_0/q) < 2^{q_0}\delta^{q_0/q}.
\end{equation*}
Setting $m = Nq_0/q$ and since we aim to minimize $q$,  we seek the largest $m$ (say $m_*$) such that
\begin{equation}\label{eq29}
    \log_2 E(q_0,m) <  q_0+\frac{m}{N}\log_2 \delta.
\end{equation}
Typical parameter choices are $q_0 = 53$ and $\delta = 10^{-16}$ which implies results in double precision, while $N$, the polynomial degree in the JFP basis required to resolve the solution to a given accuracy (which is typically the same as  $\delta$), is problem-dependent. To estimate $m_*$ in practice, we simulate $E(q_0,m)$ with $\widetilde{E}(q_0,m)$ and increase $m$ until the first value of $m$ (say $m_1$) is found for which  \eqref{eq29} is not valid, then we set $m_* = m_1-1$ and set the precision to $q = Nq_0/m_*$.

\subsection{Computational complexity}

Now we analyse the cost of the pseudo-stabilized Algorithm~2 with inputs $p$, $N$ and $\delta \ll 1$, where $p \in \mathbb{N}_+$ satisfies $\mu p = k_* \in \mathbb{N}_+$. Using \eqref{eq29} and \eqref{eq27}, it follows that $\log_2 E(1,m/q_0)<1$, so $m=\mathcal{O}(q_0)$ and since $q=\frac{q_0}{m}N$, $q=\mathcal{O}(N)$.

The first $N$ columns of $I^{(\alpha,\beta)}_{b,p,\mu}$ involve the  computation of $\mathcal{O}(N^2)$ entries, each of which requiring up to $4p+3$ multiplications, leading to a total cost of $\mathcal{O}(pN^2q\log q\log\log q)$ in $q$-bit precision\footnote{Recall that a multiplication in $q$-bit precision has $\mathcal{O}(q\log q\log\log q)$ complexity with the Schönhage--Strassen algorithm.}. Recalling that $q=\mathcal{O}(N)$, we conclude that Algorithm~2 has an overall complexity of $\mathcal{O}(pN^3\log N\log\log N)$.

\begin{remark}
	Since the accuracy drops during the recurrence, the precision can be lowered gradually to reduce the complexity by a constant factor. However, that does not change the conclusions with big O notations.
\end{remark}

\begin{remark}
    There are other computational costs that are negligible compared to the main recurrence process:
    \begin{itemize}
        \item There are divisions involving entries of the banded integration matrices which means that we can compute the inverse first with a cost of $\mathcal{O}(pN^2\log N(\log\log N)^2)$ using Newton--Raphson division.
        \item The computation of $\widetilde{E}(q_0,m)$ has $\mathcal{O}(pq_0^3\log q_0\log\log q_0)$ complexity.
    \end{itemize} 
\end{remark}

\begin{remark}
It is also possible to pseudo-stabilize Algorithm~1 (see \cref{sect:fioalgs}) since the growth of errors follow similar patterns. However, the complexity of Algorithm~1 is higher, making the overall complexity of the pseudo-stabilized version $\mathcal{O}(N^4)$. We therefore only use Algorithm~1 if $\mu$ is irrational (recall that Algorithm~2 is only applicable for rational $\mu$). 
\end{remark}